\documentclass[11pt,a4,pdftex]{article}
\usepackage{amsmath,amsthm}
\usepackage{amsfonts}
\usepackage{graphicx,float}
\usepackage[bookmarks,hyperindex,hyperindex]{hyperref}
\usepackage{fancyhdr,lastpage}
\usepackage{color}
\usepackage{lipsum}
\usepackage[framed,numbered,autolinebreaks,useliterate]{mcode}
\usepackage{bm}
\usepackage{enumerate}
\usepackage{framed}
\newlength{\leftbarwidth}
\newlength{\leftbarsep}
\hypersetup{colorlinks=true,       % false: boxed links; true: colored links
linkcolor=blue,          % color of internal links
citecolor=blue,        % color of links to bibliography
filecolor=magenta,      % color of file links
urlcolor=red  }

\newcommand{\n}{\noindent}
\newcommand{\nn}{\nonumber}

\newcommand{\drop}[1]{}

\newcommand{\im}{\mathrm{i}}

\renewcommand{\d}{\partial}
\newcommand{\rd}{\mathbb{R}^d}

\newcommand{\td}{[0, \infty)}
\renewcommand{\t}{\theta}

\newcommand{\m}{\mu}

\newcommand{\s}{\sigma}

\renewcommand{\r}{\rho}

\newcommand{\tx}[1]{\text{#1}}

\newtheorem{Def}{Definition}

\newtheorem{pro}{Proposition}
\newtheorem{Th}{Theorem}
\newtheorem{Coroll}{Corollary}
\newtheorem{Lem}{Lemma}
\newtheorem{rem}{Remark}
\newtheorem{example}{Example}
\newtheorem{Assumption}{Assumption}

\newcommand{\bD}{\begin{Def}}
\newcommand{\eD}{\end{Def}}
\newcommand{\bT}{\begin{Th}}
\newcommand{\eT}{\end{Th}}
\newcommand{\bC}{\begin{Coroll}}
\newcommand{\eC}{\end{Coroll}}
\newcommand{\beq}{\begin{equation}}
\newcommand{\eeq}{\end{equation}}
\newcommand{\bpr}{\begin{pro}}
\newcommand{\epr}{\end{pro}}
\newcommand{\fr}[2]{\frac{#1}{#2}}

\newcommand{\mc}[1]{\mathcal{ #1 }}
\newcommand{\bbra}[2]{\Big< #1, #2 \Big>}
\newcommand{\bbar}[1]{\Big\| #1 \Big\|}

\newcommand{\ee}[2]{\beq #1 #2 \eeq}
\newcommand{\ea}[2]{\begin{eqnarray} #1 #2 \end{eqnarray}}
\newcommand{\de}[2]{\begin{Def} #1 #2 \end{Def}}
\newcommand{\bP}[2]{\begin{pro} #1 #2 \end{pro}}
\newcommand{\bL}[2]{\begin{Lem} #1 #2 \end{Lem}}
\newcommand{\bTh}[2]{\begin{Th} #1 #2 \end{Th}}
\newcommand{\bRem}[2]{\begin{rem} #1 #2 \end{rem}}
\newcommand{\bA}[2]{\begin{Assumption} #1 #2 \end{Assumption}}
\newcommand{\bCo}[2]{\bC #1 #2 \eC}
\newcommand{\ex}[2]{\begin{example}#1 #2 \end{example}}
\newcommand{\pr}[1]{\begin{proof} \textrm{#1} \end{proof}}
\newcommand{\f}{\mathcal{F}}

\newcommand{\e}{\mathbb{E}}
\newcommand{\re}{\mathbb{R}}
\newcommand{\xb}[1]{\mathbf{#1}}
\newcommand{\xx}[1]{\mathbf{X}_{#1}}

\newcommand{\cov}[2]{\tx{cov}(#1, #2)}

\newcommand{\lecbbt}[1]{\lim_{n \rightarrow \infty}\mathbb{E}\Big[\sup_{t \leq T \wedge \tau_{n}^r} \Big | #1 \Big |^2\Big]}

\definecolor{light-gray}{gray}{0.85}
\definecolor{shadecolor}{rgb}{1,0.8,0.3}
\newlength{\boxwidth}
\newsavebox{\boxcontainer}

\newcommand{\nTextW}{\textwidth}

\newenvironment{Mybox}%
  {\setlength{\boxwidth}{\nTextW}
   \addtolength{\boxwidth}{-2\fboxsep}
   \addtolength{\boxwidth}{-2\fboxrule}
   \begin{lrbox}{\boxcontainer}%
     \begin{minipage}{\boxwidth}}%
  {\end{minipage}\end{lrbox}%
   \begin{center}%
     \fcolorbox{black}{light-gray}{\usebox{\boxcontainer}}
   \end{center}}

\title{Tensor Approximation of Generalized Correlated Diffusions and Functional Copula Operators}
\author{Antonio Dalessandro$^{\dagger}$\footnote{\noindent The content of this document is by no means related to any opinions that might be held or practices that might be followed at the institutions I worked for, I work for or I am associated with. I welcome comments and suggestions.}~~and Gareth W. Peters$^{\dagger}$\\
$^{\dagger}$ Department of Statistical Sciences, University College London, U.K. \\
\small{info@antonio-dalessandro.com}\\
\small{gareth.peters@ucl.ac.uk}\\
}
\date{Version 1: \today}

\begin{document}
\maketitle
\abstract{\noindent 
In this paper we derive the properties for the operator representation of the multivariate dependence theorem attributed to Sklar which describes the unique representation of general dependence structures linking marginal distributions. We develop this in the context of multidimensional correlated Markov diffusion processes, in the process defining the copula infinitesimal generator which we interpret as a functional copula specification of Sklars theorem extended to representations involving multivariate generalized diffusion processes. This allows us to develop a copula function mapping framework which we demonstrate can be accurately and efficiently obtained under a discretization scheme proposed. We show the discretized approximate copula function mapping has a limiting form which produces any desired dependence structures, induced by local correlation approximations obtained from the multivariate generalized diffusion at all local points in the state space. Hence, we achieve what we denote the functional specification of the copula mapping in operator space via tensor decomposition approaches obtained from the mimicking generalized diffusion process.

We investigate aspects of semimartingale decompositions, approximation and the martingale representation for multidimensional correlated Markov processes.  A new interpretation of the dependence among processes is given using the martingale approach. We show that it is possible to represent, in both continuous and discrete space, that a multidimensional correlated generalized diffusion is a linear combination of processes that originate from the decomposition of the starting multidimensional semimartingale. This result not only reconciles with the existing theory of diffusion approximations and decompositions, but defines the general representation of infinitesimal generators for both multidimensional generalized diffusions and as we will demonstrate also for the specification of copula density dependence structures. This new result provides immediate representation of the approximate solution for correlated stochastic differential equations. We demonstrate desirable convergence results for the proposed multidimensional semimartingales decomposition approximations.\\

\n {\bf Keywords:} Martingale Problem, Martingale Representation, Semimartingales Decomposition, Tensor algebra, Copula Functions, Copula Infinitesimal Generators, Multidimensional Semimartingales Decomposition Approximations, Diffusion Approximation Convergence.
%\tableofcontents

\newpage

%*************************************************
%*************************************************
\section{Introduction}
%*************************************************
%*************************************************
We consider the multidimensional stochastic differential equation (SDE) of the form
\ee{\label{sde}}{dX_t =  \bm{b}(t,X_t)dt + \bm{\Sigma}(t,X_t)dW_t,}
\n where $\bm{b}(t,X_t) : \td \times \rd \rightarrow \rd$ and $ \bm{\Sigma}(t,X_t): \td \times \rd \rightarrow \mathbb{R}^{d \times d} $ and we assume that $\bm{b} = (b_j)$ is continuous  vector valued function and $\bm{\Sigma} = ((\Sigma_{ij}))$ is a continuous, symmetric, nonnegative definite, $d \times d$ matrix valued function. Let the infinitesimal generator associated to the SDE of eq.(\ref{sde}) be denoted,
\ee{}{ A = \Big\{\big(f, Gf = \fr{1}{2} \sum_{ij}^d \Sigma_{ij} \d_i \d_j f + \sum_{i}^d b_i \d_i f \big): f \in  C_{c}^{\infty}(\rd) \Big\},}
\n where $C_{c}^{\infty}(\rd)$ denotes the sets of smooth functions with compact support on $\rd$.
The aim of our work is to address problems in settings which may involve potentially very high dimensional state spaces which display non-trivial dependence structures, as specified by generalized multivariate correlated diffusion processes. We propose two ways to calculate the weak solution of the SDE in eq.(\ref{sde}) in an approximate manner:
\begin{enumerate}
 \item Direct approximation of the infinitesimal generator $A$ with particular emphasis in the structure of its mixed derivative terms.
  \item Decomposition of the infinitesimal generator $A$ into orthogonal components.
\end{enumerate}
The approximation schemes we propose are based on tensor algebra decompositions such as those considered in \cite{Hackbusch2012}. The novelty we introduce consists in the introduction of new concepts like the copula infinitesimal generator, correlated tensor representation, conditional infinitesimal generator and the framework developed to perform a parametric copula function mapping as will be detailed in the remainder of the paper. In general the proposed results aim to develop a new characterization of the cross space among dimensions using tensor algebra. Furthermore, our schemes are supported by approximation and convergence results that constitute a little utilised perspective to look at solutions of multidimensional SDEs.

Our investigation is focused on aspects of the semimartingale decomposition and martingale representation for multi-dimensional correlated Markov processes. The objective is to construct a  continuous time Markov chain (CTMC) that approximates or mimicks such processes and their dependence structures induced throughout the state space, which are only implicitly defined by the joint structure of the marginal process volatility functional forms and the joint coupling of the correlation structures in the driving noise processes. Once such a mimicking process is obtained we may then transform it structure to produce a copula function mapping theorem that allows one to obtain structural characterizations of general dependence frameworks, through the specification of the mimicking multivariate diffusion process.

This work is motivated by the problems of finding the expression of a multi-dimensional CTMC that both closely follows the dynamics of the corresponding correlated Ito-processes and also can effectively deal with the representation and simulation of large dimensional processes that exhibit various correlation structures. Although the literature on Markov processes and Markov chains is very rich and mature, see \cite{RogersWilliams, EthierKurtz, KaratzasShreve}, we find that there is still room for further investigatation and characterization of multi-dimensional chains and the relationship between the correlation structure among marginal Markov chains and dependence concepts like copula functions \cite{Nelsen1999} and concordance measures of dependence \cite{McNeil2005, scarsini1984measures}. In fact these concepts have always been treated separately, and there is lack in literature of a theory that reconciles them.

Our findings and results show that our approach, based on linear and tensor algebra, is a powerful way to produce accurate solutions of multidimensional correlated SDEs that exhibit a correlation that can be fully modelled through copula functions. Specifically given a multi-dimensional Ito processes whose drift and diffusion terms are adapted processes, we show how to construct the approximated infinitesimal generator and how to characterize the process properties by its associated continuous time Markov chain (CTMC). We construct an approximated weak solution to the stochastic differential equation that weakly converges to the distribution of the multi-dimensional Ito processes.

We develop an interpretation for the correlation among processes using the martingale approach applied to the study of diffusions. The novelty is that it is possible to represent, in both continuous and discrete space, that a multidimensional correlated generalized diffusion is a linear combination of processes that originate from the decomposition of the starting multidimensional semimartingale.

The only  assumption required  by our approximation approach is that the martingale problem for the associated generator of the multidimensional Markov process is well posed. \cite{StroockVaradhan} formulated the martingale problem as a means of studying Markov processes, especially multidimensional diffusions. This approach is deemed to be  more powerful and more intrinsic than the alternative approaches represented by the Markov process approach and the Ito approach.

Our result reconciles with the existing theory of diffusion approximations and decompositions existing in the probability literature and is more closely related to the work of \cite{gyongy86} and more recently to \cite{brunick2013}. In the seminal manuscript \cite{gyongy86} considers a multi-dimensional Ito process, and constructs a weak solution to a stochastic differential equation which mimics the marginals of the original Ito process at each fixed time instant. The drift and covariance coefficients for the mimicking process can be interpreted as the expected value of the instantaneous drift and covariance of the original Ito process, conditional on its terminal value. In \cite{brunick2013} the authors extend the result of \cite{gyongy86}, proving that they can match the joint distribution at each fixed time for various functionals of the Ito process. The mimicking process takes the form of a stochastic functional differential equation and the diffusion coefficient is given by the so-called Markovian projection. In our framework we further generalize findings from \cite{brunick2013} and the mimicking process takes form of a sequence of conditional continuous time Markov chains with instantaneous drift and diffusion coefficients given by projected instantaneous local moments.

The results reported in this manuscript define the general representation of the approximated infinitesimal generators for both multidimensional generalized diffusions and for what we will define as the function copula specification.

The paper is organized as following: Section \ref{sec1} introduces the martingale problem for correlated Markov processes. In Section \ref{sec:diffApprox} we introduce and characterize the approximation schemes for the infinitesimal generator of correlated Markov processes while in Section \ref{sec3} we proposed some desirable convergence results of our approximations.

%*************************************************
%*************************************************
\section{Martingale problems for correlated Markov Processes}\label{sec1}
%*************************************************
%*************************************************
Our analysis takes place on a measure space $(\Omega, \f, \mu)$, where $\m$ is a non zero $\s$-finite, positive measure on the measurable space $(\Omega, \f)$. We denote by $\| f \|_p$ for $p \in [1, \infty]$, the $L^p(\m)$ norm of a function $f$: $\| f \|_{p}^p = \int |f|^p d\m$. Furthermore we define the Markov semigroup $(P_t(x,dy))$ on a countable set $E$ as the family of probability transition kernels on $E$ depending on the parameter $t \in [0, \infty)$, such that
\ee{}{\int_{y \in E} P_s(x, dy)P_t(y, dz) = P_{s + t}(x, dz),~~~\tx{for}~s,t \in [0, \infty),  }
and its action on bounded  and positive functions is denoted by
\ee{}{P_t f(x) = \int f(y) P_t(x, dy).}
\n The family of operators $P_t$ satisfies the following axioms $P_0 = I$, $P_t \circ P_s = P_{t + s}$, $\lim_{t \rightarrow 0^+} \| P_t f - f \| = 0$.
\n Markov processes  are naturally related to Markov semigroups because the probability measure $P_t(x, dy)$ with $P_t f(x) = \e_x[f(X_t)]$ is the law of $\{X_t\}_{t \geq 0}$ the process itself starting from the value $x$ at time $0$. We assume in this paper that the measure $\m$ is an invariant measure for the semigroup $P_t$, and this means that for all $f \in L^1(\m)$,
\[\int P_t f(x)\m(dx) = \int f(x) \mu(dx),\]
and that $P_t$ is a contraction semigroup in $L^1(\m)$ for all $t$. We define the infinitesimal generator of a strongly continuous contraction semigroup by the map
\ee{}{A f := \lim_{t \rightarrow 0^+} \fr{P_t f - f}{t},  }
\n for all $f \in D_p$, with $D_p$ a dense set of right continuous left limit (rcll) functions that is a subspace of $L^p(\mu)$.  The properties of $A$ and $D_p$ entirely specify the semigroup $P_t$. In fact given $f \in D_p$ the function $U(x,t) = P_t f(x)$ is the unique solution of the equation
\ee{\label{cauchy}}{\fr{\d U(x,t)}{\d t} = A U(x,t),}
\n defined in $D_p$ for all $t > 0$ and with $U(x,0) =  \nu$, the initial probability distribution.
\n We restrict ourself to the set of compactly supported functions denoted by $C_{c}^{\infty}(\rd)$ and belonging to the set $\bar{D} \subseteq D_p$ and we study only generators of local diffusions of the form
\ee{\label{ig}}{A_s = \sum_{i}^d b_i(s,x)\fr{\d}{\d x_i} + \fr{1}{2} \sum_{i,j}^d \Sigma_{i,j}(s,x)\fr{\d^2}{\d x_i \d x_j}},
\n where $b = (b_i(s, x))$, $i = 1, \ldots, d$ is a drift vector and $\Sigma = ((\Sigma_{i,j}))$, $i,j = 1, \ldots, d$ is a dispersion matrix with $\Sigma = C C'$ that characterize locally stochastic differential equations (SDE) with expression
\ee{\label{mainSDE}}{X_t = X_0 + \int_{0}^t b(s, X_s)ds  + \int_{0}^t C(s, X_s)dW_s,   }
\n where $W$ is a d-dimensional Brownian motion. We note that the dynamics of $X_t$ are characterized completely by the infinitesimal operator and therefore by the laws of the drift and diffusion coefficients, including the conditional probability law. In particular the infinitesimal operator in eq.(\ref{ig2}) is specified also by such coefficient with an explicit expression.
In this respect the connection between the generator of $X_t$ and the solution of the SDE for $X_t$ has a rigorous formulation given by the Martingale problem of \cite{StroockVaradhan}.
It is straightforward to show the connection between the differential operator $A$ and the probabilistic interpretation of the solution  to the corresponding SDE.
\n If $f \in C_{c}^2(\re^d)$ then Ito's Lemma yields
\ee{}{f(X(t)) = f(X(0)) + \int_{0}^t A f(X(s))ds + \int_{0}^t \nabla f(X(s))' \s(X(s)) dW(s).   }
This means that
\ee{}{f(X_t) - f(X_0) - \int_{0}^t A f(X_s)ds,}
is a local martingale. In particular the processes
\ea{}{M_t &=& X(t) - X(0) - \int_{0}^t b(X(s))ds = \int_{0}^t \s(X(s))dW(s),\\
&& M_{t}^i M_{t}^j - \int_{0}^t \s_{ij}(X_s)ds, }
are local martingales.
\de{\label{martingaleProblem}}{A probability measure $\nu$ on $(C[0, \infty)^d$, $\mc{B}(C[0, \infty)^d ))$ under which
\ee{}{M_f(t) = f(X_t) - f(X_0) - \int_{0}^t A_{s}f(X_s))ds}
is a martingale for every $f \in \bar{D} $, and is called a solution to the martingale problem associated with the operator $A_t$, where $A_t$ is as defined in eq.(\ref{ig}).}

In addition, it will also be useful to observe the following representation and uniqueness results regarding the process $X_t$ and the infinitesimal operator $A_t$, as detailed for instance in \cite{StroockVaradhan}.

\bTh{}{The process $X_t$ is a weak solution to the SDE of eq.(\ref{sde}) if and only if it satisfies the martingale problem of eq.(\ref{martingaleProblem}) with $A_t$ as the infinitesimal operator of $X_t$ as defined in eq.(\ref{ig}).  }

Another fundamental result we mention that is directly relevant background to the framework we develop is that given by the following theorem.

\bTh{}{[Uniqueness of SDE solution.]\cite[Th. 3.2.1.]{StroockVaradhan} The SDE
\ee{}{dX_t = b(X_t)dt + C(X_t)dW_t,~~~ 0 \leq t \leq s,~ X_0 = \nu,  }
with $\nu$ a random variable independent from $W_t$, $t \geq 0$ and $\e[|\nu|^2] < \infty$, has unique solution  adapted to the filtration generated by $W_t$ and $\nu$  if the measurable functions $b(x)$ and $C(x)$, for $t \leq s$, satisfy the following two conditions:
\begin{enumerate}
  \item Lipschitz continuity
  \ee{}{|b(x) - b(y)| + |C(x) - C(y)| \leq K(|x - y|),~~ \tx{for all}~x,y \in \rd, }
  \item Linear Growth
  \ee{}{|b(x)|+ |C(x)| \leq K(1 + |x|),~~ \tx{for all}~x \in \rd. }
\end{enumerate}
}
Having developed the general theoretical construct for the processes we will be working with, we may now proceed with the specification of our framework.

%*************************************************
%*************************************************
\section{Approximation of Correlated Markov Processes}\label{sec:diffApprox}
%*************************************************
%*************************************************
In this section we illustrate two new approximation schemes based on tensor algebra  for correlated Markov processes where dynamics are expressed by the SDE in eq.(\ref{mainSDE}).
The novelty we introduce in these approximations involves the characterization of the cross space among dimensions using tensor algebra and in the introduction of new concepts like the copula infinitesimal generator, correlated tensor representation, conditional infinitesimal generator and parametric copula functions mapping. The proposed schemes we develop are based on tensor algebra which makes them highly amenable to address problems in high dimesional state spaces for correlated processes. The approximations involve:
\begin{enumerate}
\item Direct approximation of the infinitesimal generator $A$ with particular emphasis on its mixed derivatives terms;
\item Decomposition of the infinitesimal generator $A$ into orthogonal components.
\end{enumerate}

In what follows we first illustrate the SDE approximation for $d = 1$, in order to establish some useful notation and the building blocks of the approximation schemes. We then present the details of approximations of correlated processes when $d \geq 2$. 

In order to approximate a one dimensional process, we construct a state space $\xb{X} \subset E$ with $n \in \mathbb{N}$ elements and define the sets of such stencils by
\ee{\label{discreteSupport}}{ \xb{X}:=\{\underline{x} := -\frac{2^{2n}}{2^n} = x_1, -\frac{2^{2n}}{2^n} + \frac{1}{2^n} := x_2, \ldots,\frac{1}{2^n},\frac{2}{2^n}, \ldots, x_m = \frac{2^{2n}}{2^n} :=   \overline{x}\} \quad \text{and} \quad \mathbf{X}^o := \mathbf{X}\backslash  \partial \mathbf{X},}
with $h = \frac{1}{2^n}$ a positive constant that represents the discretization unit, and where the boundary $\partial \mathbf{X}$ consist of the smallest (i.e. $\underline{x}$) and largest (i.e. $\overline{x}$) elements in $\mathbf{X}$, possibly a countably infinite set, and the interior $\mathbf{X}^o$ is the complement of the boundary. We denote by $\pi_n: E \rightarrow \xb{X} $ the bounded linear transformation from the continuous state space to the discretized one.

It is possible to construct the continuous time Markov chain $X^{(n)} := \{X^{n}_t\}_{t \geq 0}$ as the discrete approximation of $X$ on $\mathbf{X}$ by building a matrix $A^{(n)} = \{a(x_i,x_j)\},~~ i,j = 1, \ldots, 2^{2n +1} + 1 = m$, that is the discretized counterpart of $A$ in eq.(\ref{ig}) and each entry can be calculated by solving the following system of local moment matching equations:
\ee{\label{1dapprox}}{\left\{ \begin{array}{l}
a(x_1,x_2) = a(x_{m},x_{m-1}) = 0,\\
a(x_i,x_{i+1}) = \frac{1}{2} \Big(\frac{b(x_i)}{h} + \frac{\sigma^2(x_i)}{h^2} \Big),\\
a(x_i,x_{i-1}) = \frac{1}{2} \Big(\frac{\sigma^2(x_i)}{h^2} - \frac{b(x_i)}{h} \Big),\\
a(x_i,x_i) = - (a(x_i,x_{i-1}) + a(x_i,x_{i+1})),\\
\end{array}\right. }
for all $i = 2, \ldots, 2^{2n +1}$, with $-\frac{\sigma^2(x_i)}{h} \leq b(x_i) \leq \frac{\sigma^2(x_i)}{h}$. However, the discrete state space $\mathbf{X}$ does not need to be uniform and alternative discretization routines are presented in \cite{TavellaRandall}.

\bRem{}{In the following we will denote by $A^{(n)} = A^{(n)}_{X} := A^{(n)}_{X}(b(x_i),\s(x_i)),~ i = 1, \ldots, n,$ the approximated infinitesimal generator for the Markov process $X$ with local parameters $\m(\cdot)$ and $\s(\cdot)$. In particular $A^{(n)}_{X_k}$ is the approximated infinitesimal generator for the Markov process $X_k$. Furthermore, given a process $X_k$ we will denote by $\xb{X}_k$ the vector corresponding to its discrete state space, in a similar fashion as $\xb{X}$ in eq.(\ref{discreteSupport}). }

\bA{\label{boundaryBeh}}{
Note that for $x \in \partial \mathbf{X}$, for computational aspects, we impose an absorbing boundary condition. However, it is important to choose the boundary states sufficiently in the extreme of the state space that the laws of the processes $X^{(n)}$ and $X$ are close to each other during the finite time interval of interest in the approximation.
}

The resulting matrix $A^{(n)}$ is a tridiagonal matrix in $\mathbb{R}^{m \times m}$ , $m = 2^{2n +1} + 1 $, with always positive extra-diagonal elements.
The previous system calculates the entries of the generator by specifying the first and second instantaneous moments of the process $X^{(n)}$ that have to coincide with those of $X$ on the set $\mathbf{X}^o$. This is equivalent to satisfy the following conditions
\ee{\label{mmGen}}{\mathbb{E}_{X_t}[(X_{t + \Delta t} - X_t)^z] = \mathbb{E}_{X_t}[(X^{n}_{t + \Delta t} - X^{n}_t)^z] + o(\Delta t), \qquad z \in \{1,2\}~~\text{and}~~X^{(n)} \in  \mathbf{H}^o.}
Furthermore, one could in principle produce more accurate results by matching higher instantaneous moments of the process and in general there will be a trade off between the number of local moments matched and the coarsity of the grid/stencil $h$.

The numerical problem we face is of the same type as in eq.(\ref{cauchy}) and its analytic solution is $U_t = e^{tA}\nu$ and represents the transient probability of a Markov chain with $n$ states. 

We are now in a position to introduce the approximation schemes when $d \geq 2$. Let $A_{X_k}^{(n_k)} \in \re^{n_k \times n_k}$ denote the approximated infinitesimal generator for the continuous Markov process $k$, with $k = 1, 2, \ldots, d$. This notation is useful when describing $d$ correlated processes and the unique approximated generator for the multidimensional process. Each matrix $A_{X_k}^{(n_k)}$ is tridiagonal and its entries calculated using instantaneous local moment matching as described in in equation (\ref{1dapprox}).

The representation of the infinitesimal generator for correlated Markov processes given in eq.(\ref{ig}) can be rewritten as
\ee{\label{ig2}}{L_t = \sum_{i}^d b_i(t,x)\fr{\d}{\d x_i} + \fr{1}{2} \sum_{\substack{i,j\\ i=j}}^d a_{i,j}(t,x)\fr{\d^2}{\d x_i \d x_j} +  \fr{1}{2} \sum_{\substack{i,j\\ i \neq j}}^d a_{i,j}(t,x)\fr{\d^2}{\d x_i \d x_j}. }
\n Now denote this operator in two components,
\ee{}{A_{X_1,\ldots, X_d} := \sum_{i}^d b_i(t,x)\fr{\d}{\d x_i} + \fr{1}{2} \sum_{\substack{i,j\\ i=j}}^d a_{i,j}(t,x)\fr{\d^2}{\d x_i \d x_j},}
\n and
\ee{\label{mixedD}}{A^{(c)}_{X_1,\ldots, X_d} = \fr{1}{2} \sum_{\substack{i,j\\ i \neq j}}^d  a_{i,j}(t,x)\fr{\d^2}{\d x_i \d x_j},}
\n then  we can rewrite $L_t$ in eq.(\ref{ig2}) as the sum of two linear operators
\ee{\label{continuousOperators}}{L_t = A_{X_1,\ldots, X_d} + A^{(c)}_{X_1,\ldots, X_d}. }
In particular $A_{X_1,\ldots, X_d}$ is the continuous operator for the independent Markov processes $X_1,..., X_d$, while $A^{(c)}_{X_1,\ldots, X_d}$ is the continuous operator just for the dependence structure of such processes.

Next we develop two ways to approximate eq.(\ref{ig2}) on a discrete multidimensional space $\bigoplus_{i = 1}^n \xb{X}_i$, namely:

\begin{enumerate}
  \item With direct approximation of the operators $A^{(n_1)}_{X_1}, \ldots, A^{(n_d)}_{X_d} $ within the orthogonal dimensions and the operator $A^{(c)(n_1, \ldots, n_d)}_{X_1,\ldots, X_d}$ cross spaces. This approach directly approximates the cross derivatives operators and their action over the cross state spaces;
  \item The operator is approximated only over the orthogonal spaces. This is possible through the introduction of the notion of conditional operator $A_{X_i| X_j = a}^{(n_i)}$ with $a \in \xb{X}_j$ and $i \neq j, j = 1, \ldots, d$.
\end{enumerate}

\de{}{[Multidimensional approximated generator] The multidimensional approximated generator for $n$ independent Markov process $\{X_i(t)\}_{t\geq 0}$, for $i = 1, \ldots, d$ with approximated generators   $A^{(n_k)}_{x_k} \in \re^{n_k \times n_k}$, $k = 1, 2, \ldots, n$ is the $n_1 n_2 \cdots n_d \times n_1 n_2 \cdots n_d$ matrix
\ee{\label{multiL}}{L^{(n_1 \cdots n_d)}_{X_1, \ldots, X_d} = A^{(n_1)}_{x_1} \oplus A^{(n_2)}_{x_2} \oplus \cdots \oplus A^{(n_d)}_{x_d}.  }}

\bP{}{[Joint generator I]\label{jointRepr2}  Let $X_1$ and $X_2$ be two Markov processes with associated approximated infinitesimal generators $A^{(n_1)}_{X_1}$ and $A^{(n_2)}_{X_2}$ respectively. It is possible to define a factor $Z$, a Markov process that acts instantaneously on the spaces $(\xx{1}, \xx{2})$ with  generator $A_{Z}^{(n_1 n_2)}$,  such that the infinitesimal approximated generator of the correlated processes $(X_1, X_2)$ with local correlation parameter $\r := \r(\xx{1}, \xx{2})$ can be written as
\ee{\label{cGen}}{A^{(n_1 n_2)}_{X_1,X_2} = A^{(n_1)}_{X_1} \oplus A^{(n_2)}_{X_2} + A^{(c)(n_1 n_2)}_{X_1, X_2}}
\n where
\ea{\label{copulaOperator}}{ A^{(c)(n_1 n_2)}_{X_1, X_2} & = & - \Big(L^{(n_1)}_{Z} \oplus L^{(n_2)}_{Z} \Big) - \tx{diag} \Big(I^{(n_1)} \otimes diag(L^{(n_2)}_{Z})\Big) \nn \\
 &+& 1_{\{ \r > 0\}} \Big( I^{(n_2)}_m \otimes L^{(m)(n_1)}_{Z} +  L^{(p)(n_2)}_{Z} \otimes I^{(n_1)}_p \Big)\nn \\
 &+& 1_{\{ \r < 0\}} \Big( I^{(n_2)}_m \oplus L^{(p)(n_1)}_{Z} +  L^{(p)(n_2)}_{Z} \otimes I^{(n_1)}_p \Big).}
Furthermore eq.(\ref{cGen}) can be rewritten differently instead as a conditional structure given by
\ee{\label{cGen2}}{A^{(n_1 n_2)}_{X_1,X_2} = A^{(n_1)}_{X_1|Z} \oplus A^{(n_2)}_{X_2|Z} + A^{(n_1 n_2)}_{Z}. }}
\pr{First we consider the mixed derivative terms in eq.(\ref{mixedD}) when $d = 2$. For positive correlation the following approximations hold for finite element approximations for the partial derivative operators,
\ea{\label{localMixed}}{\fr{\d f}{\d x_1 \d x_2} &=& \fr{\Big(\fr{\d f}{\d x_2}\Big)_{i+1,j} - \Big(\fr{\d f}{\d x_2}\Big)_{i,j}  }{\Delta x_1} = \fr{f_{i+1,j+1} - f_{i+1,j} - f_{i,j+1} + f_{i,j}}{\Delta x_1 \Delta x_2} + O(\Delta x_{1}) + O(\Delta x_{2}),\\
\fr{\d f}{\d x_1 \d x_2} &=& \fr{\Big(\fr{\d f}{\d x_2}\Big)_{i,j-1} - \Big(\fr{\d f}{\d x_2}\Big)_{i-1,j-1}  }{\Delta x_1} = \fr{f_{i,j} - f_{i,j-1} - f_{i-1,j} + f_{i-1,j-1}}{\Delta x_1 \Delta x_2} + O(\Delta x_{1}) + O(\Delta x_{2}).\nn }
\n Equations (\ref{localMixed}) can be combined to yield,
\ee{\label{mixedP}}{\fr{\d f_{ij}}{\d x_1 \d x_2} = \fr{f_{i+1,j+1} - f_{i+1, j} - \big(f_{i ,j-1} - 2 f_{i,j} + f_{i, j+1} \big) - f_{i-1, j} + f_{i-1, j-1}}{2 \Delta x_1 \Delta x_2} + O(\Delta x_{1}^2) + O(\Delta x_{2}^2) + O(\Delta x_1  \Delta x_2) }
Note the same scheme applies for negative correlation.
However eq.(\ref{mixedP}) can be decomposed into the following three terms $T3 - (T1 + T2)$, where
\ea{}{T1 &=& \fr{f_{i+1,j} - 2 f_{i,j} + f_{i-1,j} }{2 \Delta x_1 \Delta x_2}, \nn \\
T2 &=& \fr{f_{i,j+1} - 2 f_{i,j} + f_{i,j-1} }{2 \Delta x_1 \Delta x_2},\nn \\
T3 &=& \fr{f_{i+1,j+1} - 2 f_{i,j} + f_{i-1,j-1} }{2 \Delta x_1 \Delta x_2}.}
Next we observe that the operator $A^{(c),(n_1 n_2)}_{X_1, X_2}$ is a \lq correlation' operator acting on the joint discretized product space $\xx{1} \times \xx{2}$. The term T1 acts only along the discretized support of $X_1$, the term T2 acts only along the discretized support of $X_2$, while T3 acts only along the discretized cross support for both $X_1$ and $X_2$. We then use these finite difference operators to calculate the entries of the operators in eq.(\ref{copulaOperator}). In particular, we use the scheme T1 for  $L^{(n_1)}_{Z}$, the scheme T2 for $L^{(n_2)}_{Z}$  and T3 for $1_{\{ \r > 0\}} \Big( I^{(n_2)}_m \otimes L^{(m)(n_1)}_{Z} +  L^{(p)(n_2)}_{Z} \otimes I^{(n_1)}_p \Big)$. The magnitude of the local instantaneous intensities is $\r(\xx{1},\xx{2}) \s{\xx{1}}, \s{\xx{2}}$.
\n We can therefore rewrite eq.(\ref{cGen}) as
\ea{}{A^{(n_1 n_2)}_{X_1,X_2} &=& (A^{(n_1)}_{X_1} - L^{(n_1)}_{Z}) \oplus (A^{(n_2)}_{X_2} - L^{(n_2)}_{Z}) \nn\\
&-& \Big[\tx{diag}(I^{(n_1)} \otimes diag(L^{(n_2)}_{Z})) \nn \\
 &+& 1_{\{ \r > 0\}} \Big( I^{(n_2)}_m \otimes L^{(m)(n_1)}_{Z} +  L^{(p)(n_2)}_{Z} \otimes I^{(n_1)}_p \Big)\nn \\
 &+& 1_{\{ \r < 0\}} \Big( I^{(n_2)}_m \oplus L^{(p)(n_1)}_{Z} +  L^{(p)(n_2)}_{Z} \otimes I^{(n_1)}_p \Big)\Big]. }
If we define the following operators,
\ea{}{A^{(n_1)}_{X_1|Z} &=& (A^{(n_1)}_{X_1} - L^{(n_1)}_{Z}), \nn \\
A^{(n_2)}_{X_2|Z} &=& (A^{(n_2)}_{X_2} - L^{(n_2)}_{Z}),  \nn \\
A^{(n_1 n_2)}_{Z} &=& - \tx{diag}(I^{(n_1)} \otimes diag(L^{(n_2)}_{Z})) + 1_{\{ \r > 0\}} \Big( I^{(n_2)}_m \otimes L^{(m)(n_1)}_{Z} +  L^{(p)(n_2)}_{Z} \otimes I^{(n_1)}_p \Big)\nn \\
&& + 1_{\{ \r < 0\}} \Big( I^{(n_2)}_m \oplus L^{(p)(n_1)}_{Z} +  L^{(p)(n_2)}_{Z} \otimes I^{(n_1)}_p \Big).}This proves eq.(\ref{cGen2}).}

\bP{}{[Multivariate Uniform Distribution Infinitesimal Operator]\label{UniformIG} The operator $A^{(c)(n_1 n_2)}_{X_1, X_2}$ is the infinitesimal generator associated to a bivariate correlated distribution with uniform marginals.}
\pr{Let's consider the operator $A^{(c)(n_1 n_2)}_{X_1, X_2}$ in eq.(\ref{copulaOperator}) acting on the joint Hilbert space $H_{X_1 X_2} = H_{X_1} \oplus H_{X_2}$. Without loss of generality we consider the structure of the approximated infinitesimal operator when $\r >0$, and we note that the case with $\r < 0$ is identical. In this case the operator is given by,
\ee{\label{copulaG_poisson}}{A^{(c)(n_1 n_2)}_{X_1, X_2} = - \Big(L^{(n_1)}_{Z} \oplus L^{(n_2)}_{Z} \Big) - \tx{diag} \Big(I^{(n_1)} \otimes diag(L^{(n_2)}_{Z})\Big) +  \Big( I^{(n_2)}_m \otimes L^{(m)(n_1)}_{Z} +  L^{(p)(n_2)}_{Z} \otimes I^{(n_1)}_p \Big),}
\n and it is a linear function of the operators $L^{(n_1)}_Z$ and $L^{(n_2)}_Z$ that act on the spaces $H_{X_1}$ and $H_{X_2}$ respectively.  }

\bP{}{[Conditional Infinitesimal Generator]\label{jointRepr1} Let the infinitesimal generators of generalized diffusions $\{X_i(t)\}_{t \geq 0}, ~ i = 1, 2, \ldots$ be
\[ A_{X_i} f(x) = \mu_i(x) \frac{\partial f(x)}{\partial x} + \fr{1}{2} \sigma_{i}^2(x)  \frac{\partial^2 f(x)}{\partial x^2} , ~~\tx{i} = 1, 2, \ldots\] and with general and local properties as described in section \ref{sec:diffApprox}. Let's assume that the diffusions $\{X_i(t)\}_{t \geq 0}, ~ i = 1, 2, \ldots$ are locally correlated with instantaneous local covariation given by
\ee{}{\bbra{\s_i(x) dW_i(t)}{\s_j(y) dW_j(t)} = \s_i(x)\s_j(x)\r_{ij}(x,y), ~~~~~ \tx{for all}~i,j = 1, 2, \ldots}
In order to introduce the concept of conditional infinitesimal generator, without loss of generality, we consider the 2-dimensional operator,
\ee{\label{continuousOperators2}}{L_t = A_{X_1, X_2} + A^{(c)}_{X_1, X_2} }
being the d-dimensional case just an algebraic extension. We will make all the necessary dimensionality considerations for the d-dimensional representation.
Therefore let's consider the processes pair $X_1 = X$ and $X_2 = Y$. The conditional approximated infinitesimal generator $A^{(n_1)}_{X | Y}$ is defined by the sequence of operator matrices $\Big\{A^{(n_1)}_{X | Y = y_j} \in \re^{n_1 \times n_1} \Big\}, ~~ y_j \in \xb{Y}$ each of whose entries are obtained according to local moment matching by:
\ee{\label{conditionalIG}}{A^{(n_1 n_2)}_{X_1 | X_2 = y_j} = \left\{ \begin{array}{l}
a(x_1,x_2) = a(x_{m},x_{m-1}) = 0,\\
a(x_i,x_{i+1}) = \frac{1}{2} \Big(\frac{\mu_1(x_i) + \r_{12}(x_i, y_j) \fr{\s_1(x_i)}{\s_2(y_j)}(y_j - \m_2(y_j))}{h} + \frac{\sigma_{1}^2(x_i)(1 - \r_{12}^2(x_i, y_j) ) }{h^2} \Big),\\
a(x_i,x_{i-1}) = \frac{1}{2} \Big(  \frac{\sigma_{1}^2(x_i)(1 - \r_{12}^2(x_i, y_j) ) }{h^2}  - \frac{\mu_1(x_i) + \r_{12}(x_i, y_j) \fr{\s_1(x_i)}{\s_2(y_j)}(y_j - \m_2(y_j))}{h} \Big),\\
a(x_i,x_i) = - (a(x_i,x_{i-1}) + a(x_i,x_{i+1})),\\
\end{array}\right. }
for all $y_j \in \xb{Y}$, $x_i \in \xb{X}$, with $-\frac{\sigma_{k}^2(x_i)}{h} \leq \mu_k(x_i) \leq \frac{\sigma_{k}^2(x_i)}{h}$ for $k = 1,2$, $ -1 \leq \r_{12}(x_i, y_j) \leq 1$, and $a(x_i,x_j) \geq 0$ for all $i \neq j$.}
\pr{Let's consider the Markov processes pair $(X,Y)$ and how to derive its corresponding infinitesimal generator approximation in matrix form. The local instantaneous intensities are calculated in the same way as described in Section \ref{sec:diffApprox} for the one dimensional case, namely using local moment matching as reported in eq.(\ref{mmGen}). For the 2-d case the instantaneous intensities are calculated using the local transition  kernel
\ee{\label{2dPr}}{p(X_{t + \Delta t}, Y_{t +  \Delta t}| X_{t}, Y_{t}),~~\tx{as}~ \Delta t \rightarrow 0.}
In particular considering the local states $x_i \in \xb{X}$ and $y_j \in \xb{Y}$ the transition probability in eq.(\ref{2dPr}) can be rewritten in local form as
\ee{}{p(x_{i+1},y_{j+1}| x_{i}, y_j) = p(x_{i+1}| y_{j+1}, x_{i}) p(y_{j+1}| y_j), }
\n where
\ee{}{\Big(Y_{t^+} = y_{j+1}| Y_{t} = y_{j}\Big ) \sim  p(y_{j+1}| y_j) = N \Big(\m(y_j), \s^2(y_j)\Big)}
\ea{\label{instCondM}}{\Big(X_{t^+} = x_{i+1} | Y_{t} = y_{j + 1}, X_t = x_i\Big ) &\sim&  p(x_{i+1}| y_{j+1}, x_{i})\\ &=& N\Big(\m(x_i) + \fr{\s(x_i)}{\s(y_j)}\r(x_i,y_j)(y_j - \m(y_j)), (1 - \r^2(x_i,y_j) \s^2(x_i)) \Big)\nn}
\n The calculation of the matrix entries for $A^{n_1 n_2}_{X_1 | X_2 = y_j}$ is done in identical way as described in Section \ref{sec:diffApprox}, imposing the instantaneous moments of eq.(\ref{instCondM}) for all $y_j \in \xb{Y}$ and this completes our proof. Extension to dimensions larger than two is straightforward due to independence of each conditional operator.}

\n From proposition \ref{jointRepr1} it is clear that we can represent the approximated multidimensional generator of eq.(\ref{multiL}) as a decomposition of independent conditional generators. Furthermore it is also possible to exploit standard results of conditional probability partitioning in order to facilitate the local characterization of the independent multidimensional conditional generators. Given a multivariate gaussian variable $\bm{X} \sim N(\bm{\m}, \bm{\Sigma})$  and consider the partition of $\bm{X}$ and equivalently of $\bm{\m}$ and $\bm{\Sigma}$ into
\ee{}{\bm{X} = \begin{bmatrix}
  \bm{x_1} \\
  \bm{x_2}
 \end{bmatrix}, ~~ \bm{\m} = \begin{bmatrix}
  \bm{\m_1} \\
  \bm{\m_2}
 \end{bmatrix}, ~~ \bm{\Sigma} = \begin{bmatrix}
  \bm{\Sigma_{1,1}} & \bm{\Sigma_{1,2}} \\
  \bm{\Sigma_{2,1}} & \bm{\Sigma_{2,2}}
 \end{bmatrix}   }
\n Then $(\bm{x_1}| \bm{x_2}) =  \bm{y} = \bm{x_1} + \bm{C} \bm{x_2} $, where $\bm{C} = - \bm{\Sigma_{1,2}}\bm{\Sigma_{2,2}}^{-1}$, the conditional distribution of the first partition given the second, is $N(\overline{\bm{\m}}, \overline{\bm{\Sigma}})$, with mean
\ee{}{\overline{\bm{\m}} =  \bm{\m_1} + \bm{\Sigma_{1,2}}\bm{\Sigma_{2,2}}^{-1} (\bm{x_2} - \bm{\m_2}), }
and covariance matrix
\ee{}{\overline{\bm{\Sigma}} = \bm{\Sigma_{1,1}} + \bm{\Sigma_{1,2}}\bm{\Sigma_{2,2}}^{-1}\bm{\Sigma_{2,1}}.}

More generally if we denote by $\Sigma$ the covariance matrix introduced within the definition of the operator in eq.(\ref{ig}), and by $\bm{X} \sim N(\bm{\m}, \bm{\Sigma})$ a multivariate normal vector, and fixed $t \geq 0$ the covariance matrix $\Sigma = \{a_{i,j}(t,x)\} \in \re^{d \times d} $  of equations (\ref{ig}), (\ref{ig2}) is
\ee{}{\bm{\Sigma} =
 \begin{pmatrix}
  \Sigma_{1,1} & \Sigma_{1,2} & \cdots & \Sigma_{1,d} \\
  \Sigma_{2,1} & \Sigma_{2,2} & \cdots & \Sigma_{2,d} \\
  \vdots  & \vdots  & \ddots & \vdots  \\
  \Sigma_{d,1} & \Sigma_{d,2} & \cdots & \Sigma_{d,d}
 \end{pmatrix}.
}

Conditional probability partitioning is a very important property when creating the sequence of conditional approximated generators as in eq.(\ref{conditionalIG}) because large multivariate Gaussian vectors can be easily partitioned as the combination of sets of independent sub-multivariate Gaussian vectors. Each of the sub-multivariate Gaussian vectors can be further characterized and locally approximated through a principal component analysis (PCA). Therefore it is possible to construct a reduced dimensionality infinitesimal generator of a large dimension process without an aggregate PCA of the global process covariance structure.

\ex{}{[3-D approximated generator]\label{ex_3d} Given the results in propositions (\ref{jointRepr1}) we show how to calculate $A^{(n_1 n_2 n_3)}_{X_1 X_2 X_3}$. In fact we can express the approximated 3-D generator under a conditional decomposition according to
\ee{}{A^{(n_1 n_2 n_3)}_{X_1 X_2 X_3} = A^{(n_1)}_{X_1 | X_2 X_3} \oplus A^{(n_2)}_{X_2 | X_3} \oplus A^{(n_3)}_{X_3}. }
Note that $A^{(n_1)}_{X_1 | X_2 X_3} = A^{(n_1)}_{X_1 | X_2 = x_2, X_3 = x_3} \in \re^{n_1 \times n_1} $ with $x_2 \in \bm{X_2}$ and  $x_3 \in \bm{X_3}$.
  }

\bP{}{[Joint Transition Kernel I] The transition kernel solution of the Cauchy problem in eq.(\ref{cauchy}) for the infinitesimal generator in eq.(\ref{cGen}) is given by
\ee{\label{JTKI}}{P^{(n_1 n_2)}(\xx{1},\xx{2}) = P^{(n_1 n_2)}(\xx{1},\xx{2}; \r = 0)P^{(c)(n_1 n_2)}(\xx{1},\xx{2}; \r) = \sum_{\xb{Z}}P(\xx{1}|\xb{Z}) \otimes P(\xx{1}|\xb{Z}) P(\xb{Z}) } }
\pr{Due to the linearity of the operator in eq.(\ref{cGen2})
\ee{}{A^{(n_1 n_2)}_{X_1,X_2} = A^{(n_1)}_{X_1|Z} \oplus A^{(n_2)}_{X_2|Z} + A^{(n_1 n_2)}_{Z},}
\n we obtain
\ea{}{e^{A^{(n_1 n_2)}_{X_1,X_2}} &=& e^{A^{(n_1)}_{X_1|Z} \oplus A^{(n_2)}_{X_2|Z} + A^{(n_1 n_2)}_{Z},}\nn\\
&=& \sum_{\xb{Z}} e^{A^{(n_1)}_{X_1|Z} \oplus A^{(n_2)}_{X_2|Z}} e^{A^{(n_1 n_2)}_{Z}} }
which proves the rhs of eq.(\ref{JTKI}). Analogously the linearity of the operator in eq.(\ref{cGen}) implies
\ea{}{e^{A^{(n_1 n_2)}_{X_1,X_2}} &=& e^{ A^{(n_1)}_{X_1} \oplus A^{(n_2)}_{X_2} + A^{(c)(n_1 n_2)}_{X_1, X_2}}\nn\\
&=& e^{A^{(n_1)}_{X_1} \oplus A^{(n_2)}_{X_2}} e^{A^{(c)(n_1 n_2)}_{X_1, X_2}}, }
which proves the lhs of eq. (\ref{JTKI}). Due to the equivalence of the operator representation in eq.(\ref{cGen}) and eq.(\ref{cGen2}), this proves equivalence of the transition probability kernels in eq.(\ref{JTKI}).}

\bL{}{[Equivalence of the joint representations]  Given two correlated Markov processes $X_1$ and $X_2$ the approximated joint transition probability kernel $P^{(n_1 n_2)}_{X_1,X_2}(\xx{1}, \xx{2})$ is the solution of the Cauchy problem of eq.(\ref{cauchy}) where the infinitesimal operator can be expressed in either cross space and marginals decomposition or as conditional decomposition given by
\ee{\label{g1}}{A^{(n_1 n_2)}_{X_1,X_2} = A^{(n_1)}_{X_1} \oplus A^{(n_2)}_{X_2} + A^{(c)(n_1 n_2)}_{X_1, X_2},}
\n or
\ee{\label{g2}}{A^{(n_1 n_2)}_{X_1,X_2} = A^{(n_1)}_{X_1 | X_2} \oplus A^{(n_2)}_{X_2}. }
\n The approximations in eq.(\ref{g1}) and eq.(\ref{g2}) are equivalent.       }

\bTh{}{[Correlated kernel tensor representation]\label{GeneralSDEsolution} Given $m$ correlated Markov processes the approximated solution of the martingale problem of eq.(\ref{martingaleProblem}) is given by the product measure representation involving the generator decomposition given according to,
\ee{\label{generalSolution}}{P^{(n_1 \cdots n_m)}_{X_1, \ldots, X_m}(t) = e^{t A^{(n_1)}_{X_1 |X_2, \ldots, X_m}} \otimes e^{t A^{(n_2)}_{X_2 |X_3, \ldots, X_m}} \otimes \cdots \otimes e^{t A^{(n_m)}_{X_m}}} }

\pr{In order to calculate the m-dimensional solution of eq.(\ref{generalSolution}) we exploit the orthogonality of the conditional approximated infinitesimal operators $A^{(n_i)}_{X_i, \ldots, X_m}$ for all $i$.
This follows from propositions \ref{jointRepr1}. Then the transition density of eq.(\ref{generalSolution}) is computed as illustrated in example  \ref{ex_3d}.  }

It is important at this stage to introduce the explicit expressions of the quadratic variation and quatrain covariation among CTMC with generators approximated using the method introduced in proposition (\ref{jointRepr1}).

\bL{}{[Orthogonal CTMCs] \label{orthogonalCTMC} Let $\overline{\bm{X}} = [ X^{(n_1)}_1, \ldots, X^{(n_d)}_d ]'$ be a d-dimensional CTMC and define its partition as $\overline{\bm{X}} = [\bm{x_1}, \bm{x_2}]'$. Let's denote define $\bm{y} = \bm{x_1} + \bm{C} \bm{x_2} $, where $\bm{C} = - \bm{\Sigma_{1,2}}\bm{\Sigma_{2,2}}^{-1}$. Then the chains $\bm{x_2}$  and the conditional chains $\bm{y}$ are orthogonal and $\tx{var}(\bm{x_1} | \bm{x_2}) = \tx{var}(\bm{y}) $. }
\pr{We want to show that the instantaneous covariation of the chains $\bm{x_2}$ and $\bm{y}$ is zero. Therefore we compute
\ea{}{\cov{\bm{x_2}}{ \bm{y}} &=& \cov{\bm{x_2}}{ \bm{x_1}} + \cov{\bm{x_2}}{ \bm{C x_2}} \nn\\
&=& \bm{\Sigma_{1,2}} + \bm{C} \cov{\bm{x_2}}{ \bm{x_2}} \nn \\
&=& \bm{\Sigma_{1,2}} - \bm{\Sigma_{1,2}}\bm{\Sigma_{2,2}}^{-1}\bm{\Sigma_{2,2}} = 0  }
\n Furthermore we have
\ea{}{\tx{var}(\bm{x_1} | \bm{x_2}) &=& \tx{var}(\bm{x_1} + \bm{C} \bm{x_2}) \nn \\
&=& \tx{var}(\bm{x_1}) + \bm{C} \tx{var}(\bm{x_2})\bm{C'} + \bm{C}\cov{\bm{x_1}}{ \bm{x_2}} + \cov{\bm{x_2}}{ \bm{x_1}}\bm{C'} \nn \\
&=& \tx{var}(\bm{y} | \bm{x_2}) = \tx{var}(\bm{y}).}
 }

\n Given the results from proposition  \ref{jointRepr1} and Lemma \ref{orthogonalCTMC} with some algebra it is straightforward to obtain the covariance matrix $\bm{\Sigma}$. We have
\ee{}{[\bm{y}, \bm{x_2}] \begin{bmatrix}
  \overline{\bm{\Sigma}} & 0 \\
  0 & \bm{\Sigma_{2,2}}
 \end{bmatrix} \begin{bmatrix}
  \bm{y} \\
  \bm{x_2}
 \end{bmatrix} = [\bm{x_1}, \bm{x_2}] \begin{bmatrix}
  \bm{\Sigma_{1,1}} & \bm{\Sigma_{1,2}} \\
  \bm{\Sigma_{2,1}} & \bm{\Sigma_{2,2}}
 \end{bmatrix}   \begin{bmatrix}
  \bm{x_1} \\
  \bm{x_2}
 \end{bmatrix}.
}

%*************************************************
\subsection{Functional Copula Constructions}
%*************************************************
We are interested in generating joint distributions with a variety of dependence structures and to achieve this purpose we use copula function specifications. In what follows we explore how copulas and families of copulas are generated and introduce a general approach to construct copulas in a tensor product space. In particular we show how copulas are related to other methods of generating joint distributions in a Hilbert space based on specified marginal tensors and decomposition properties of the product space. In doing this we exploit some desirable convergence results of the proposed transformed tensor representation of multivariate correlated processes to the continuous copula functions.

After the definition of a copula function we recall Sklar's theorem, a fundamental result about the relationship between marginals and joint distribution for multivariate correlated random variables. This is important to recall as we will develop a new characterization of Sklar theorem as a function generator representation.

\de{}{[m-Copula Function]\label{copulaFunction} An m-dimensional copula (or m-copula) is a function $C$ from the unit m-cube $[0,1]^m$ to the unit interval $[0,1]$ which satisfies the following conditions:
\begin{enumerate}
  \item $C(1,\ldots, 1 ,a_n, 1,\ldots, 1) = a_n$ for every $n \leq m$ all $a_n$ in $[0,1]$;
  \item $C(a_1,\ldots,a_m) = 0$ if $a_n = 0$ for $a_n \leq a_m$;
  \item $C$ is $m$-increasing.
\end{enumerate}}

 \bRem{}{Property 1 says that if the realizations of $m-1$ variables are known each with marginal probability one, then the joint probability of the outcomes is the same as the probability of the remaining uncertain outcome. Property 2 is sometimes referred to as the grounded property of a copula. It says that the joint probability of all outcomes is zero if the marginal probability of any outcome is zero. Property 3 says that the C-volume of any $m$-dimensional interval is non-negative. Properties 2 and 3 are general properties of multivariate cdfs. It follows that an $m$-copula can be defined as an $m$-dimensional cdf whose support is contained in $[0,1]^m$ and whose one-dimensional margins are uniform on $[0,1]$. In other words, an $m$-copula is an $m$-dimensional distribution function with all $m$ univariate margins being $U(0,1)$.}

The relationship between distribution functions and copulas is given by the following result, see \cite{Sklar...}.

\bTh{}{[Sklar's Theorem] Let $X$ and $Y$ be random variables with distribution functions $F$ and $G$ respectively and joint distribution function $H$. Then there exists a copula $C$ such that for all
$(x,y) \in \re \times \re$
\ee{\label{sklarDistrF}}{H(x,y) = C(F(x),G(y))}
If $F$ and $G$ are continuous, then $C$ is unique; otherwise, $C$ is uniquely determined on \hbox{Ran(F) $\times$ Ran(G)}. Conversely, if $C$ is a copula and $F$ and $G$ are distribution functions, then the function $H$ defined by (\ref{sklarDistrF}) is a joint distribution function with margins $F$ and $G$.}

By Sklar's theorem, given continuous margins $F_1$ and $F_2$ and the joint continuous distribution function $F(x_1,x_2)$ = $C(F_1(x_1), F_2(y_2))$, the corresponding copula is generated using the unique inverse transformation
\ee{}{C(u_1, u_2) = C(F(x_1), F_2(x_2)) = F(x_1, x_2) = F(F_{1}^{-1}(u_1),F_{2}^{-1}(u_2)),}
where $u_1$ and $u_2$ are standard uniform variates.
\n Given the result in theorem \ref{GeneralSDEsolution} it is possible to extend it to the representation of joint distribution function using tensor algebra. Without loss of generality we present the result in two dimensions. Extension to higher dimensional case is straightforward by induction.

\bP{}{[Sklar's Theorem in Generator Space]\label{SklarGenerator} Let $P^{(n_1 n_2)}_{X_1, X_2}(\xx{1}, \xx{2})(t), t \geq 0$, be the approximated transition probability kernel solution of the martingale problem in eq.(\ref{martingaleProblem}) for the infinitesimal  generator
\ee{}{L(s) = \sum_{i}^2 b_i(s,x)\fr{\d}{\d x_i} + \fr{1}{2} \sum_{i,j}^2 \Sigma_{i,j}(s,x)\fr{\d^2}{\d x_i \d x_j}.}
\n Then the approximated  joint distribution  function for the process starting at $(\underline{x}_1, \underline{x}_2)$ is given by
\ee{\label{approxF}}{F^{(n_1 n_2)}_{X_1,X_2}(\xx{1},\xx{2}| X_1 = \underline{x}_1, X_2 = \underline{x}_2 ) = F^{(n_1)}_{X_1}(\xx{1} |X_1 = \underline{x}_1) \otimes F^{(n_2)}_{X_2}(\xx{2} |X_2 = \underline{x}_2) } }
\pr{\ea{}{P^{(n_1 n_2)}_{X_1, X_2}(\xx{1}, \xx{2}|X_1 = \underline{x}_1, X_2 = \underline{x}_2 )(t) &=& \Big(e^{A^{(n_1)}_{X_1}t} \otimes e^{A^{(n_2)}_{X_2}t } \Big) \mathbf{1}_{\{X_1 = \underline{x}_1, X_2 = \underline{x}_2 \}} \nn \\
 &=& e^{A^{(n_1)}_{x_1}t}\mathbf{1}_{\{X_1 = \underline{x}_1\}} \otimes e^{A^{(n_2)}_{x_2}t}\mathbf{1}_{\{X_2 = \underline{x}_2\}}\nn \\
  &=&  P^{(n_1)}_{X_1}(\xx{1} | X_1 = \underline{x}_1) \otimes P^{(n_2)}_{X_2}(\xx{2} | X_2 = \underline{x}_2) }
\n due to the conditional independence of the marginal tensors. Then
\ee{}{ \sum_{i = \underline{x}_1 }^{x_1} P(x^{(1)}_i | X_1 = \underline{x}_1) \otimes \sum_{j = \underline{x}_2 }^{x_2}P(x^{(2)}_j | X_2 = \underline{x}_2) = F^{(n_1 n_2)}_{X_1,X_2}(\xx{1},\xx{2}| X_1 = \underline{x}_1, X_2 = \underline{x}_2 ). } }

Due to the general validity of our result we can proceed to formulate it in the following theorem.

\bTh{}{[Joint Distribution Function Convergence]\label{copulaConvergence} Let $X_1(t)$ and $X_2(t)$ be correlated Markov processes with marginal distributions $F_{X_1}$ and $F_{X_2}$ and approximated joint distribution function $F^{(n_1 n_2)}_{X_1,X_2}(\xx{1}, \xx{2})$ as in eq.(\ref{approxF}).
Let $C(u_1, u_2) : [0, 1]^2 \mapsto [0, 1]$ a continuous copula function. Then the following convergence result holds
\ee{}{\lim_{n_1,n_2 \rightarrow \infty} F^{(n_1 n_2)}_{X_1,X_2}(\xx{1}, \xx{2}) = C(F_{X_1}(X_1), F_{X_2}(X_2)) } }

\bCo{}{[General Copula Mapping] Let $D^{(n_1, n_2, \ldots)} := [\fr{i}{n_1}] \times [\fr{j}{n_2}] \times \ldots $ where $i = 0, \ldots, n_1 $ and $j = 0, \ldots, n_2$ and so on, denote the discretization of the unit hypercube. Let $C^{n_1,n_2}(u_i,v_j)$ be a copula distribution function defined for all $(u_i, v_j) \in D^{(n_1, n_2)}$. Then according to theorem \ref{copulaConvergence}, the following equalities hold:
\ee{}{C^{(n_1,n_2)}(u_i,v_j) = C^{(n_1,n_2)}\Big(F_{X_1}^{(n_1)}(x^{(1)}_i),F_{X_2}^{(n_2)}(x^{(2)}_j)\Big ) = F(x^{(1)}_i, x^{(2)}_j)}
\n Furthermore, if we denote by $\theta$ the set of the copula parameters, and by $\textbf{p}$ the set of local cross space parameters for the approximated tensor representation, then given  $C^{(n_1,n_2)}\Big(F_{X_1}^{(n_1)}(x^{(1)}_i),F_{X_2}^{(n_2)}(x^{(2)}_j); \theta \Big )$, its equivalent representation in a tensor space is given by
\ee{\label{cop_map}}{ \min_\textbf{p}\bbar{C^{(n_1,n_2)}(u_i,v_j; \t) -  F^{(n_1 n_2)}_{X_1,X_2}(\xx{1},\xx{2}; \textbf{p})}_{2}^2 }
}
\n Eq.(\ref{cop_map}) means that for a given parametric copula distribution function $C$ belonging to any copula family it is possible to find a set of local parameters $\textbf{p}$ that would produce a joint distribution function $F$ that has minimal Euclidean distance from $C$ in the tensor space.
\n A practical way to solve eq.(\ref{cop_map}) is to compute the local likelihood with respect to the set of parameters $\textbf{p} = \{p_{ij}\}$ such that
\ee{\label{logLikelihood}}{\min_{p_{ij}} \Big \{\log\Big(C_{\tx{target}}(u_i, u_j; \t) - 
F^{(n_1 n_2)}_{X_1,X_2}\left(F_{X_1}^{(n_1)}(x^{(1)}_i),F_{X_2}^{(n_2)}(x^{(2)}_j); p_{ij}\right)
\Big \}, ~~\tx{for all}~i,j. }
\n where $C_{\tx{target}}$ can be any target copula function.

\bTh{}{[Copula Infinitesimal Operator]\label{copulaGenerator}
Let $c^{(n_1,n_2)}(u_i,v_j)$ be a copula density function defined for all $(u_i, v_j) \in D^{(n_1, n_2)}$. Given the result in theorem \ref{copulaConvergence}, we have that
\ee{}{c^{(n_1,n_2)} \Big(F_{X_1}^{(n_1)}(x^{(1)}_i),F_{X_2}^{(n_2)}(x^{(2)}_j) \Big) = P^{(c)}(x^{(1)}_i, x^{(2)}_j),~~ \tx{for all}~i,j. }
\n with infinitesimal operator given by eq.(\ref{copulaOperator}).}

\bTh{}{[Copula Tensor Representation] Given the tensor product basis $Z = B_X \oplus B_Y$ associated to $A^{(n_x,n_y)}_{X,Y}$ as in eq.(\ref{g1}), where $B_X$ and $B_Y$ denote the basis of the operators $A^{(n_x)}_{X}$ and $A^{(n_x)}_{X}$ respectively, it is possible to specify a point $(x_i, y_j) \in (\xb{X}, \xb{Y})$ with corresponding associated subspace $(M \oplus M^\perp) \subset Z$ of the direct sum of vectors $x \in M$ and $y \in M^\perp$ with origin $(x_i, y_j)$, such that the following relations hold,
\ee{\label{direcsum-marginals}}{z = x \oplus y}
\n and
\ee{\label{direcsum-marginals-corr}}{z = z' + z'' =  (x \oplus y) + (-a \oplus -a^{\perp}),}
where the vector $a$ is the shared component among the vectors $x$ and $y$ or the instantaneous local covariance part of the joint process $(X_t, Y_t)$ and represents the copula.  Eq.(\ref{direcsum-marginals}) and eq.(\ref{direcsum-marginals-corr}) are represented in fig.(\ref{TensorDecomposition}).

\begin{figure}[htbp]
  \centering
  \includegraphics[scale=.8]{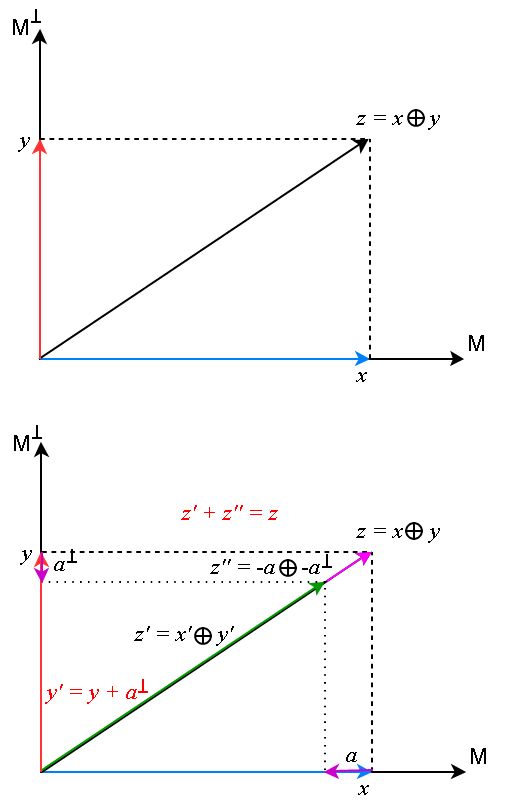}\\
  \caption{Representation of $z = z' + z'' =  (x \oplus y) + (-a \oplus -a^{\perp})$.}\label{TensorDecomposition}
\end{figure}
}

%*************************************************
%*************************************************
\section{Convergence of the Approximated Generator}\label{sec3}
%*************************************************
%*************************************************
In section \ref{sec:diffApprox} we defined the Markov chain $X_{t}^{(n)}$ approximating the multidimensional generalized diffusion $X_t$. Here we present weak convergence results for $X_{t}^{(n)}$ to the solution of the SDE for $X_t$, introduced in eq.(\ref{sde}). The convergence is studied from different perspectives: first from a semigroup point of view, secondly in a spectral way through a Fourier unitary transformation of the approximated generator, and lastly through the martingale problem for the associated infinitesimal generator.

We want that the continuous Markov chain $X^{(n)} := \{X^{n}_t\}_{t \geq 0}$ has a dynamics as close as possible to the corresponding approximated process $\{X_t\}_{t \geq 0}$. At this purpose we can define an error
\ee{\label{errGenerator}}{ \epsilon_n(f) := \sup_{x \in \mathbf{H}} \| A^{(n)} f(x) - Af(x) \|.}
Using the semigroup approach to weak convergence we can state that, if $\epsilon_n(f)$ tends to zero as $n$ tends to infinity for $f$ in $\overline{\mathcal{D}}$, then the sequences of processes $X^{(n)}$, converges weakly to $X$, in the space $D_\mathbf{E}$. The space $D_\mathbf{E}$, is the space of right continuous functions $f: \mathbf{E} \rightarrow \mathbb{R}$ with left limits, see Ch. 5, pp. 115 of \cite{EthierKurtz} for more details. 
\n The following theorem states that if the error $\epsilon_n(f)$ goes to zero as $n$ tends to infinity, implying norm convergence of the approximated generator $A^{(n)} f(x) $ to $A f(x)$, this would imply also convergence of the corresponding approximated semigroup $P^{(n)}_{t}$ to $P_t f$ and of the chain $X^{(n)}$ to $X$ for $t \geq 0$.

\bTh{\label{Th1}}{ Let $X$ be a Feller process with state space $\mathbf{E}$ and infinitesimal generator $A$ satisfying the properties and assumptions defined in section \ref{sec1}, and $X^{(n)}$ be a sequence of Markov chains with generator matrices $ A^{(n)}$. Then it holds that
\begin{equation}
\lim_{n \rightarrow \infty} \epsilon_n(f) = 0
\end{equation}
for every function $f$ in the core of $A$, which is equivalent to the statement that
\begin{equation*}
P^{(n)}_{t} \pi_n f \rightarrow P_t f, \qquad \text{for all}~ f \in C_0(\mathbf{E}),~ t \geq 0.
\end{equation*}
Then the sequences of processes $X^{(n)}$, converges weakly to $X$, in the space $D_\mathbf{E}$.
}

\pr{From Theorem 6.1 in Chapter 1 pag. 28 of \cite{EthierKurtz}, we have that if $P^{(n)}_{t}, ~ n = 1, 2, \ldots$, and $P_{t}$ are strongly continuous contraction semigroups on $B_n$ and $B$ with generators $A^{(n)}$ and $A$ respectively, and let $\overline{\mathcal{D}}$ be a core for $A$, then the following are equivalent:
\begin{itemize}
\item[(i)] For each $f \in B$
\begin{equation*}
P^{(n)}_{t} \pi_n f \rightarrow P_t f
\end{equation*}
uniformly on bounded intervals.
\item[(ii)] For each $f \in B$
\begin{equation*}
P^{(n)}_{t} \pi_n f \rightarrow P_t f
\end{equation*}
for all $t \geq 0$.
\item[(iii)] For each $f \in \overline{\mathcal{D}}$, there exists $f_n \in D(A^{(n)})$  ($D(A^{(n)})$ is the domain of $A^{(n)}$ ) for each $n \geq 1$ such that $f_n \rightarrow f$ and $A^{(n)} f_n \rightarrow Af$.
\end{itemize}
\noindent Furthermore, following \cite{EthierKurtz}, Chapter 4 pag. 172 Theorem 2.11, we have that if $P_t$ is a Feller Semigroup on $C_0(\mathbf{E})$ and that for each $t \geq 0$ and $f \in C_0(\mathbf{E})$
\begin{equation*}
P^{(n)}_{t} \pi_n f  \rightarrow P_t f.
\end{equation*}
If $X^{(n)}_0$ has limiting distribution $\nu$, then there is a Markov Process $X$ corresponding  to $P_t$ with initial distribution $\nu$ and sample paths in $D_p$, and
\begin{equation*}
X^{(n)} \rightarrow X \quad \text{in}~D_p.
\end{equation*}
}

\n Therefore, in order to be able to use the result of theorem \ref{Th1}, we just need to prove the convergence of the proposed approximated generator $A^{(n)}$ and this is done in the following theorem.

\bTh{\label{convergenceAn}}{ For all $f \in \overline{\mathcal{D}}$, with $\overline{\mathcal{D}}$ core of the generator $A$ as previously defined,
\begin{equation*}
\lim_{n \rightarrow \infty}\sup_{x \in \mathbf{H}}|A^{(n)}f(x) - Af(x)| = 0
\end{equation*}
This mode of convergence is called strong convergence.
}
\pr{ Let $f \in \overline{\mathcal{D}}$ and $x \in \mathbf{H}$, and $h := \Delta x \in  \mathbb{R}$ the discretization unit.
We can locally write the generator of the semigroup $P^{(n)}_{t}$ as
\begin{equation*}
A^{(n)} f(x) = a^{(n)}(x, x + h) \Big(f(x + h) - f(x)\Big)  + a^{(n)}(x, x - h) \Big(f(x - h) - f(x)\Big)
\end{equation*}
By Taylor approximation
\begin{equation*}
f_n(x) \approx f(x + \Delta x) \approx  f(x) + f'(x)\Delta x + \frac{1}{2}f''(x)(\Delta x)^2 + o((\Delta x)^2)
\end{equation*}
we obtain for all $i = 1, 2, \ldots$, $x_i \in \mathbf{H}$,
\begin{eqnarray*}
A^{(n)} f_n(x_i) &=&  a(x_i,x_{i+1}) \Big( f(x_i) + f'(x_i)\Delta x + \frac{1}{2}f''(x_i)(\Delta x)^2 + o((\Delta x)^2) - f(x_i) \Big)\\
 &+&  a(x_i,x_{i - 1}) \Big( f(x_i) - f'(x_i)\Delta x + \frac{1}{2}f''(x_i)(\Delta x)^2 + o((\Delta x)^2) - f(x_i) \Big)\\
&=& f'(x_i)\Delta x \Big(a(x_i,x_{i+1}) - a(x_i,x_{i-1})\Big) + \frac{1}{2}f''(x_i)(\Delta x)^2 \Big(a(x_i,x_{i+1}) - a(x_i,x_{i-1})\Big)\\
&+&  \Big(a(x_i,x_{i+1}) - a(x_i,x_{i-1})\Big)o((\Delta x)^2).
\end{eqnarray*}
Due to the fact that $f \in \overline{\mathcal{D}}$, the error term $o((\Delta x)^2)$ is uniform in $x$.
We have that
\begin{eqnarray*}
a(x_i,x_{i+1}) + a(x_i,x_{i-1}) &=& \frac{\sigma^2(x_i)}{h^2},\\
a(x_i,x_{i+1}) - a(x_i,x_{i-1})&=& \frac{\mu(x_i)}{h}.\\
\end{eqnarray*}
In section \ref{sec:diffApprox} we made precise assumptions, see \ref{boundaryBeh} on the behavior at the boundary of the process. Without loss of generality we can assume that the boundary is at a point of infinity, i.e. not attainable in finite time and it is furthermore absorbing.
It is possible to express $A^{(n)} f_n(x)$ for all $x \in \mathbf{H}$, as
\begin{eqnarray*}
A^{(n)} f_n(x) &=& f'(x)\Delta x \Big(\frac{\mu_{i}}{\Delta x}\Big) + \frac{1}{2}f''(x)(\Delta x)^2 \Big(\frac{\sigma^2(x_i)}{(\Delta x)^2}\Big) + \Big(\frac{\mu(x_i)}{\Delta x}\Big)o((\Delta x)^2)\\
&=& A f(x) + \Big(\frac{\mu(x_i)}{\Delta x}\Big)o((\Delta x)^2).
\end{eqnarray*}
We obtain
\beq
\sup_{x \in \mathbf{H}}|A_n f(x) - A f(x)| = \sup_{x \in \mathbf{H}} |\frac{\mu(x_i)}{\Delta x} o((\Delta x)^2)| = C_1 o((\Delta x)^2)  \xrightarrow{n \rightarrow \infty} 0.
\eeq
In case the elements of $\partial \mathbf{H}$, consisting of the smallest (i.e. $\underline{x}$) and largest (i.e. $\overline{x}$) elements in $\mathbf{H}$, are absorbing states of the Markov chain we can continue the above analysis with a further investigation of the weak convergence. Furthermore the behavior on the boundary of the functions $f \in \overline{D}$ is $f'(\underline{x}) = f'(\overline{x}) = 0$.
For $x_i = \underline{x}$ we have,
\begin{eqnarray*}
|A^{(n)} f(\underline{x})  - A f(\underline{x})| &\leq& |A^{(n)} f(\underline{x})  - A^{(n)} f(\underline{x} + h)|\\
&+& |A^{(n)} f(\underline{x}+ h)  - A f(\underline{x} + h)|\\
&+& |A f(\underline{x}+ h)  - A f(\underline{x})|
\end{eqnarray*}
The second term tends to 0 as shown above, the third term by continuity of
$A f$ in $\underline{x}$. For the first term we have
\begin{eqnarray*}
|A^{(n)} f(\underline{x})  - A^{(n)} f(\underline{x} + h)| &=& |a(\underline{x}, \underline{x} + h) (f(\underline{x} + h) - f(\underline{x}))\\
&-& a(\underline{x} + h, \underline{x}) (f(\underline{x} + 2h) - f(\underline{x} + h))\\
&-& a(\underline{x} + h, \underline{x} + 2h) (f(\underline{x} + h) - f(\underline{x}))|\\
&=& O(n)o(h) \xrightarrow{n \rightarrow \infty} 0.
\end{eqnarray*}
since $a(\underline{x}, \underline{x} + h)$, $a(\underline{x} + h, \underline{x})$, $a(\underline{x} + h, \underline{x} + 2h)$ are on order $n$ and $(f(\underline{x} + h) - f(\underline{x}))$, $(f(\underline{x} + 2h) - f(\underline{x} + h))$, $(f(\underline{x} + h) - f(\underline{x}))$ are of order $o(h)$ because $f'(\underline{x}) = 0$ for $f \in \overline{D}$. The result for $\overline{x}$ follows in the same way.}

%*************************************************
\subsection{Correlated Diffusions Approximation Convergence}
%*************************************************
For this purpose we need to introduce multiplication operators that are considered as an infinite-dimensional generalization of diagonal matrices and they are extremely simple to construct. Furthermore, they appear naturally in the context of the Fourier transform or when one applies the spectral theorem and deals with spectral representation of operators on Hilbert spaces. It is an equivalent way to represent the same operator and can be useful for calculations and further analysis.

We recall that the definitions of discrete Fourier transform matrix and its inverse as a unitary operator are given by,
\begin{eqnarray*}\label{UnitaryOP}
F_{s,k} &=& \frac{1}{\sqrt{n}} e^{- \im \frac{2 \pi }{n} s k },\\
F^{-1}_{k,s} &=& \frac{1}{\sqrt{n}} e^{\im \frac{2 \pi}{n} s k},\\
F_{l,k}F^{-1}_{k,j} &=& \delta_{l,j},
\end{eqnarray*}
where $\delta_{l,j}$ is the Kronecker delta, so these matrices give the resolution of the identity matrix and define a unitary transformation. Also, if $f(x)$ is a function belonging to the space $L^2$ and $B_n$ be the Brillouin zone defined as $B_n = \Big\{- \frac{\pi}{h} + kh, \quad k = 0, \ldots, \frac{2\pi}{h^2} = \frac{n}{h} \Big\}$  the transformation $\mathcal{F}_n: L^2(\mathbf{H}) \mapsto L^2(B_n)$
\begin{equation}
\mathcal{F}_n(f)(s) = \sum_{x \in \mathbf{H}}F_{s,k} f(x),
\end{equation}
is the discrete Fourier transform of $f$, with  $\mathbf{H} := h\mathbb{Z}$ \footnote{$\mathbf{H}$ is possibly unbounded but in all our practical applications we consider $\mathbf{H} \subset D := [-K, K] \subset \mathbb{R}$, $K \in (0,\infty)$, with $D$ the operator domain and assuming for simplicity also periodic boundary conditions.} .
In fact,
\begin{eqnarray}
\mathcal{F}(f)(s)&=& \frac{1}{\sqrt{n}} \int_{-\frac{\pi}{h}}^{\frac{\pi}{h}}e^{- \im s x} f(x) dx\\
&\approx& \frac{1}{\sqrt{n}} \sum_{k=0}^{n-1} e^{\im s \big(-\frac{\pi}{h} + kh\big)} \underbrace{f(-\frac{\pi}{h} + kh)}_{f_k},
\end{eqnarray}
where $h = \frac{2\pi}{n}$. We now extend the above transformation to the d-dimensional case, and the following results set the notation for our subsequent theorems and proofs.

\bL{}{For any nonnegative-definite symmetric matrix $\bm{\Lambda}$ the function
 \ee{}{\varphi_{\xb{X}} = \exp \Big( - \im \bm{ \m s} - \fr{1}{2} \bm{s' \Lambda s} \Big)t  }
\n is the characteristic function at time $t > 0$ of the random vector $\xb{X}$ with $\e[\xb{X}] = \bm{\m}$ and $\tx{Cov}[\xb{X}] = \xb{\Lambda}$.}
\bTh{}{(cf.\cite{Lukacs1958}, Th. 3.2.3) \label{jumpDistrib} Let $f(s)$ be an arbitrary characteristic function. For every real $x$ the limit
\begin{equation}
p(x) = \lim_{T \downarrow 0} \frac{1}{2T} \int_{-T}^T e^{- \im s x} f(s) dx
\end{equation}
exists and is equal to the saltus of the distribution function of $f(s)$ at the point $x$.}

We can obtain a spectral representation of the operator $L^{n_1 \cdots n_d}_{X_1, \ldots, X_d}$ by applying the above unitary transformation leading to the following diagonal operator,
\begin{equation}\label{SpectralOP}
q_{n_1 \cdots n_d}(\bm{s}) =  \mathcal{F}(L^{(n_1 \cdots n_d)}_{X_1, \ldots, X_d}(f)) \mathcal{F}^{-1}(\bm{s},\bm{s})
\end{equation}
\n where the approximated operator of eq.(\ref{multiL}) can be written as
\ee{}{L^{(n_1 \cdots n_d)}_{X_1, \ldots, X_d} = \bm{\m' \nabla} + \fr{1}{2} \bm{\s' \tilde{H} \s}}
\n with $\bm{\nabla}$ the discrete d-dimensional operator gradient and $\bm{\tilde{H}}$ the discrete d-dimensional Hessian operator.

In order to derive some converge properties of the operator (\ref{multiL}) we do this by comparing the spectral representations of the probability density functions for the continuous infinitesimal generator and its approximated counterpart, namely
\begin{equation*}
P_{t}^{(n_1 \cdots n_d)}(\bm{x},\bm{y}) =  \frac{1}{(2 \pi)^d} \int_{[-\frac{\pi}{h}, \frac{\pi}{h}]^d}e^{q_{n_1 \cdots n_d}(\bm{s})t}e^{\im \bm{s} (\bm{y} - \bm{x})}d\bm{s},
\end{equation*}
and
\begin{equation*}
p_t(\bm{x},\bm{y}) =  \frac{1}{(2 \pi)^d} \int_{\rd} \varphi_{\xb{X}} e^{\im \bm{s} (\bm{y} - \bm{x})}d\bm{s}.
\end{equation*}
In this way we are able to assess the order of convergence of the error
\begin{equation}\label{convergenceErr}
\bm{\epsilon_n} := \Big| p_t(\bm{x},\bm{y}) - P_{t}^{(n_1 \cdots n_d)}(\bm{x},\bm{y})\Big|.
\end{equation}

To assess the rate of convergence of  eq.(\ref{convergenceErr}), we exploit the relationship between the distribution function and its corresponding characteristic function and in particular we refer to the Continuity Theorem.

\bTh{}{[Continuity Theorem,  cf.\cite{Lukacs1958}, Th. 3.6.1.]\label{levyContinuityT} Let $\{F_n(x)\}$ be a sequence of distribution functions and denote by $\{f_n(s)\}$ the sequence of the corresponding characteristic functions. The sequence $\{F_n(x)\}$ converges weakly to a distribution function $F(x)$ if, and only if, the sequence $\{f_n(s)\}$ converges for every $s$  to a function $f(s)$ which is continuous at $s = 0$. The limiting function is then the characteristic function of $F(x)$.}

We can therefore focus on the analysis of the passage to the limit $\bm{h} \rightarrow 0$
of the following spectral representation, conditional on a time $t$,
\begin{equation}\label{AA}
\lim_{\bm{h} \downarrow 0}  \frac{1}{(2 \pi)^d} \int_{[-\frac{\pi}{h}, \frac{\pi}{h}]^d}e^{q_{n_1 \cdots n_d}(\bm{s})t}e^{\im \bm{s} (\bm{y} - \bm{x})}d\bm{s} =  \frac{1}{(2 \pi)^d} \int_{\rd} \varphi_{\xb{X}} e^{\im \bm{s} (\bm{y} - \bm{x})}d\bm{s}.
\end{equation}

\bTh{}{[Convergence of the d-dimensional approximated operator]For all $\bm{x}$ we consider the sequence of distribution functions $P_{t}^{(n_1 \cdots n_d)}(\bm{x},\bm{y})$  and by $\{ e^{q_{n_1 \cdots n_d}(\bm{s})t} \}$ the sequence of the corresponding characteristic functions. The sequence $\{ F_{\bm{n}}(\bm{x})\}$ converges weakly to a distribution function $F(\bm{x})$ if, and only if, the sequence $\{ e^{q_{n_1 \cdots n_d}(\bm{s})t} \}$ converges for every $\bm{s}$  to a function $\varphi_{\xb{X}}$ which is continuous at $\bm{s} = 0$. The limiting function is then the characteristic function of $F(\bm{x})$.}
\begin{proof}
We prove convergence and characterization of the rate of convergence for $d = 1$ and $d = 2$ being the proof in higher dimensions just an algebraic extension of the case of $d = 2$.  The calculation of $q_{n_1 \cdots n_d}(\bm{s})t$ is straightforward and it is just an application of the shift theorem. For $d = 1$
\begin{equation*}
\mathcal{F}_n(A^{(n_1)}_{X_1})(f)(s) =  \sum_{x \in \mathbf{H}} F_{s,k} \Big(\mu \nabla_{h_1}(f)(x) +  \frac{\sigma^2}{2} \triangle_{h_1} (f)(x)\Big)
\end{equation*}
\begin{eqnarray*}
\mathcal{F}_n(\mu \nabla_{h_1})(f)(s) &=&  \sum_{x \in \mathbf{H}} F_{s,k} \mu \nabla_{h_1}(f)(x)\\
&=& \mu \sum_{k =0}^{n-1} F_{s,k} \frac{f(kh + h_1) - f(kh - h_1)}{2h_1}\\
&=& \frac{\mu}{2h_1} \Big( e^{\im h_1 s} \sum_{k =0}^{n-1} F_{s,k} f(kh)  - e^{-\im h_1 s}\sum_{k =0}^{n-1} F_{s,k} f(kh) \Big)\\
&=& \mu \frac{e^{\im h_1 s} - e^{- \im h_1 s}}{2h_1} \mathcal{F}(f)(s) = -\im \mu \frac{\sin h_1 s}{h_1}\mathcal{F}(f)(s).
\end{eqnarray*}
Doing a similar calculation for $\mathcal{F}(\frac{\sigma^2}{2} \triangle_{h_1})(f)(s)$ we obtain
\begin{equation}\label{SpectralOP}
q_{n_1}(s) =  \Big( -\im \mu \frac{\sin h_1 s}{h_1} + \sigma^{2} \frac{\cos(h_1 s) -1}{h_{1}^2}\Big) (s).
\end{equation}
We have
\begin{eqnarray*}\label{ContinuousKerH}
P_{t}^{(n_1)}(x,y) &=& p^{(n_1)}(y \leq X_{t} \leq  y + h_1| X_0 = x)\\
&=& \frac{1}{n_1} \sum_{s \in B_n}e^{q_{n_1}(s)t}e^{\im s (y-x)}.
\end{eqnarray*}
Without loss of generality, this result is a particular case of the continuity theorem (\ref{levyContinuityT}), and the convergence error of eq.(\ref{convergenceErr}) is measured in correspondence of saltus point of the distribution, see Th.(\ref{jumpDistrib}).
Lets consider the integral on the left end side of eq.(\ref{AA}) for $d = 1$ and it can be rewritten as
\begin{eqnarray*}
 \frac{1}{2 \pi} \int_{-\frac{\pi}{h_1}}^{\frac{\pi}{h_1}}e^{\Big( -\im \mu \frac{\sin h_1 s}{h_1} + \sigma^{2} \frac{\cos(h_1 s) -1}{h_{1}^2}\Big) (s)t}e^{\im s (y-x)}ds &=&   \frac{1}{2 \pi} \int_{-\frac{\pi}{h_1}}^{-K}e^{\Big( -\im \mu \frac{\sin h_1 s}{h_1} + \sigma^{2} \frac{\cos(h_1 s) -1}{h_{1}^2}\Big) (s)t}e^{\im s (y-x)}ds\\
 && + \frac{1}{2 \pi} \int_{-K}^{K}e^{\Big( -\im \mu \frac{\sin h_1 s}{h_1} + \sigma^{2} \frac{\cos(h_1 s) -1}{h_{1}^2}\Big) (s)t}e^{\im s (y-x)}ds\\
   && + \frac{1}{2 \pi} \int_{K}^{\frac{\pi}{h_1}}e^{\Big( -\im \mu \frac{\sin h_1 s}{h_1} + \sigma^{2} \frac{\cos(h_1 s) - 1}{h_{1}^2}\Big) (s)t}e^{\im s (y-x)}ds
\end{eqnarray*}
\noindent and it is possible to make the first and the third integral on the right end side arbitrary small by choosing a large number $K > 0$ and by selecting $h_1 >0$ sufficiently small. If we consider the second integral we can analyze the behaviour as $h_1 \rightarrow 0$.
We notice that the function
\begin{eqnarray*}
\lim_{h_1 \rightarrow 0} \Big( \im \mu  \frac{\sin h_1 s}{h_1}\Big) &=& \lim_{h_1 \rightarrow 0}  \im \mu~\frac{s}{s} \frac{\sin h_1 s}{h_1}\\
&=& \lim_{h_1 \rightarrow 0}  \im  s \mu  \frac{\sin h_1 s}{h_1 s}= \im  s \mu
\end{eqnarray*}
and
\begin{eqnarray*}
\lim_{h_1 \rightarrow 0} \Big( \sigma^2  \frac{\cos h_1 s - 1}{h_{1}^2}\Big) &=& \lim_{h_1 \rightarrow 0}  \sigma^2 \frac{s^2}{s^2} \frac{\cos h_1 s - 1}{h_{1}^2}\\
&=& \lim_{h_1 \rightarrow 0}  \sigma^2 s^2 \frac{\cos h_1 s - 1}{(h_1 s)^2}= -\frac{1}{2}\sigma^2 s^2.
\end{eqnarray*}
\noindent We would like to examine in more details the order of convergence of the above functions as $h_1 \rightarrow 0$. For the limit
\[\lim_{h_1 \rightarrow 0}  \frac{\sin h_1}{h_1} = 1, \]
using $\sin h_1 = h_1 - \frac{h_{1}^3}{6} + \ldots$, we get
\[ \frac{\sin h_1}{h_1} - 1  = \frac{\sin h_1 - h_1}{h_1} = -\frac{h_{1}^3}{6h} + \ldots = -\frac{{h_1}^2}{6} + \ldots\]
we find an order of $O(h_{1}^2)$. In the same way using $\cos h_1 = 1 - \frac{h_{1}^2}{2} + \ldots$, we can assess the order of convergence of the limit
\[\lim_{h_1 \rightarrow 0}  \frac{1 - \cos h_1}{h_{1}^2} = \frac{1}{2}, \]
\[ \frac{1 - \cos h}{h_{1}^2} = \frac{1 - \Big( 1 - \frac{h_{1}^2}{2} + \ldots \Big)}{h_{1}^2} = \frac{1}{2} + \ldots \]
and convergence order of $O(1)$.
Therefore the order of convergence is at most $O(h_{1}^2)$. This results can be extended to all the marginals of a d-dimensional approximated operator in case of independent marginals.
In the presence of correlation  we have the presence of mixed derivative terms. For $d = 2$ the calculation of $q_{n_1 n_2}$ is as follows:
\ea{\label{2dF}}{
\f_{\bm{n}} \Big(A^{(n_1 n_2)}_{X_1, X_2}\Big)(f)(\bm{s}) &=& \f_{\bm{n}} \Big( A^{(n_1)}_{X_1} \oplus A^{(n_2)}_{X_2} + A^{(c)(n_1 n_2)}_{X_1, X_2} \Big)(f)(\bm{s}) \\
&=& \sum_{x_1 \in \xb{X}_1} \sum_{x_2 \in \xb{X}_2} F_{s,k} \Big(\big( \mu_1 \nabla_{h_1} +  \frac{\sigma_{2}^2}{2} \triangle_{h_1} + \mu_2 \nabla_{h_2} +  \frac{\sigma_{2}^2}{2} \triangle_{h_2} +  \r \sigma_{1} \sigma_{2} \nabla_{h_1}\nabla_{h_2} \big) (f)(\bm{x})\Big). \nn
}
Eq. (\ref{2dF}) is equivalent to
\ee{}{\f_{\bm{n}} \Big( A^{(n_1)}_{X_1} \oplus A^{(n_2)}_{X_2} + A^{(c)(n_1 n_2)}_{X_1, X_2} \Big)(f)(\bm{s}) = \f_{n} \Big( A^{(n_1)}_{X_1}\Big)(f)(\bm{s}) \oplus \f_{n} \Big(A^{(n_2)}_{X_2}\Big)(f)(\bm{s}) + \f_{\bm{n}} \Big(A^{(c)(n_1 n_2)}_{X_1, X_2} \Big)(f)(\bm{s}). \nn}
Therefore it is sufficient to analyze the term
\ee{}{\f_{\bm{n}} \Big(A^{(c)(n_1 n_2)}_{X_1, X_2} \Big)(f)(\bm{s}),}
where the operator $A^{(c)(n_1 n_2)}_{X_1, X_2}$ was introduced in  proposition \ref{jointRepr2} as the correlation operator.

\ea{}{\f_{\bm{n}} \Big(A^{(c)(n_1 n_2)}_{X_1, X_2} \Big)(f)(\bm{s}) &=& \sum_{x_1 \in \xb{X}_1} \sum_{x_2 \in \xb{X}_2} F_{s,k} \Big( \r \sigma_{1} \sigma_{2} \nabla_{h_1}\nabla_{h_2} \big) (f)(\bm{x}) \Big) \nn \\
&=& \r_{12} \sigma_{1} \sigma_{2} \fr{1}{h_1 h_2}\Big( \cos(h_1 s_1 + h_2 s_2)  - \cos(h_1 s_1)  - \cos(h_2 s_2) + 1   \Big),
}
because of the mixed derivative term approximation of proposition \ref{jointRepr2}. The following approximations hold:
\ea{}{\cos(h_1 s_1 + h_2 s_2) - 1 &=& - \fr{h^{2}_2 s^{2}_2}{2} - s_1 s_2 h_1 h_2 - \fr{h^{2}_1 s^{2}_1}{2} + \fr{h^{2}_1 s^{2}_1 h^{2}_2 s^{2}_2}{4} + ( 1 + h_1 + h_{1}^2) O(h_2)^3 + O(h)^3, \nn \\
\cos(h_1 s_1) - 1 &=& - \fr{s_{1}^2 h_{1}^2}{2} + O(h_1)^3, \nn \\
\cos(h_2 s_2) - 1 &=& - \fr{s_{2}^2 h_{2}^2}{2} + O(h_2)^3.}
Therefore,
\ee{}{\lim_{h_1, h_2 \rightarrow 0} \r_{12} \sigma_{1} \sigma_{2} \fr{1}{h_1 h_2}\Big( \cos(h_1 s_1 + h_2 s_2)  - \cos(h_1 s_1)  - \cos(h_2 s_2) + 1   \Big) = - \r_{12} \sigma_{1} \sigma_{2} s_1 s_2,}
that is exactly the characteristic function of the covariance term for the bivariate normal distribution at time $t$. Given this proof, extension to higher dimensions is algebraically straightforward.
\end{proof}

%*************************************************
\subsection{Weak Convergence of the Approximating Chain Using the Martingale Central Limit Theorem}
%*************************************************
Here we present a weak convergence result for the multivariate diffusion approximation introduced in Section \ref{sec:diffApprox}, along the lines of the general diffusion convergence Theorem 4.1 at pag. 354 of \cite{EthierKurtz}. The convergence result we propose is based on the arguments belonging to the formulation of diffusion theory in terms of the martingale problem, see \cite{StroockVaradhan}, which requires minimal assumptions about the smoothness of the coefficients of the SDE and can be seen as an extension of the martingale central limit theorem [Martingale CLT, Th. 1.4, \cite{EthierKurtz}, pag. 339]\label{martingaleCLT}. The results are mainly obtained by compactness arguments which do not require a priori regularity. These arguments are the same as those devised to provide an existence theory for the multidimensional SDE under consideration, and they just refer to properties of the SDE coefficients, that are the prerequisites for our analysis as per section \ref{sec1}.

\bTh{}{ Let $\bm{\Sigma} = ((\Sigma_{ij}))$ be a continuous, symmetric,nonnegative definite, $d \times d$ valued function on $\rd$ and let $b: \rd \rightarrow \rd$ be continuous. Let
\ee{}{ A = \Big\{\big(f, Gf = \fr{1}{2} \sum \Sigma_{ij} \d_i \d_j f + \sum b_i \d_i f \big): f \in  C_{c}^{\infty}(\rd) \Big\}}
and suppose that the $C_{\rd}[0, \infty)$ martingale problem for A is well posed.
\n For $n = 1, 2, \ldots$ let $X_n$ and $B_n$ be processes with sample paths in $D_{\rd}[0, \infty)$ and let $A_n = ((A_{n}^{ij}))$ be a symmetric $d \times d$ matrix-valued  process such that $A_{n}^{ij}$ has sample paths in $D_R[0, \infty)$ and $A_n(t) - A_n(s)$ is non negative definite for $t > s \geq 0$. Set $\f_{t}^n = \s\Big( X_n(s), B_n(s), A_n(s): s \leq t \Big)$.
\n Let $\tau_{n}^r = \inf \{t: |X_n(t)| \geq r~~or~~|X_n(t-)| \geq r \}$, and suppose that
\ee{\label{m1}}{M_n = X_n - B_n}
and
\ee{}{M^{i}_n M^{j}_n - A^{ij}_n}
are $\f^{n}_t$-local martingales, and for each $r>0$, $T>0$, and $i,j = 1,2, \ldots, d$
\ee{\label{m2}}{\lecbbt{X_n(t) - X_n(t-)} = 0,}
\ee{\label{m3}}{\lecbbt{B_n(t) - B_n(t-)} = 0,}
\ee{}{\lecbbt{A^{ij}_n(t) - A^{ij}_n(t-)} = 0,}
\ee{}{\sup_{t \leq T \wedge \t_{n}^r} \Big | B^{i}_n(t) - \int_{0}^t b_i(X_n(s))ds \Big |  \xrightarrow{p} 0 }
\n and
\ee{}{\sup_{t \leq T \wedge \t_{n}^r} \Big | A^{ij}_n(t) - \int_{0}^t \Sigma_{ij}(X_n(s))ds \Big | \xrightarrow{p} 0.}
\n Suppose that $P(X^{(n)}(0)^{-1}) \Rightarrow \nu$, with $\nu$ the starting distribution. Then $X^{(n)} \Rightarrow X$ , with $X$ the solution of the martingale problem for $(A, \nu)$.   }

\pr{ Due to the fact that the martingale problem for $A$ is well posed, the process $M(t) = X(t) - B(t)$ is a martingale. By the optional stopping theorem, see \cite{RogersWilliams}, if $\tau$ is a stopping time also $M(\tau)$ is a martingale. In particular this is valid also for the stopping time $\tau_{n}^r = \inf \{t: |X_n(t)| \geq r~~or~~|X_n(t-)| \geq r \}$. By eq.(\ref{m1}) the process $M_n(t) = X_n(t) - B_n(t)$ is a martingale. Relatively compactness properties of $M_n(t)$ implies relatively compactness of $X_n(t)$ and $B_n(t)$ and therefore for $X_n(t \wedge \tau_n)$ and $B_n(t \wedge \tau_n)$ as in eq.(\ref{m2}) and eq.(\ref{m3}) respectively.
This extends also to the stopped martingale $M_n(t \wedge \tau_n \wedge \tau_a )$ where $\tau_a = \inf\{t : A_n^{ii}(t) > t \sup_{|x| \leq r} a_{ii}(x) + 1, \tx{for~some~i} \}$.
Furthermore relative compactness in a set $C$ is a condition equivalent to the condition that each sequence in $C$ contains a convergent subsequence, see for example \cite{Billingsley}.
\n This means that every subsequence $X_{n_k}(t \wedge \tau^{r}_n)  \Rightarrow X^{\bar{r}}(t \wedge \tau^r$), where $\tau^r = \inf\{t : |X^{\bar{r}}(t)| \geq r \}$, for all $\bar{r} \geq r > 0$
Therefore the stopped processes
\ee{}{M^{\bar{r}}(t \wedge \tau^r) = X^{\bar{r}}(t \wedge \tau^r) - \int_{0}^{t \wedge \tau^r} b(X^{\bar{r}}(s))ds }
\ee{}{M_{i}^{\bar{r}}(t \wedge \tau^r)M_{j}^{\bar{r}}(t \wedge \tau^r)  - \int_{0}^{t \wedge \tau^r} a_{ij}(X^{\bar{r}}(s))ds }
are martingales and by Ito's lemma
\ee{}{f(X^{\bar{r}}(t \wedge \tau^r))  - \int_{0}^{t \wedge \tau^r} A f(X^{\bar{r}}(s))ds }
\n is a martingale for each $f \in C^{\infty}_c(\rd)$, and $A f(X^{\bar{r}}(s))$ is the approximated infinitesimal generator applied  to the function $f$.
\n In particular if the martingale problem is well posed uniqueness argument for the solution hold and hold also for the solution for the stopped problem, hence
\ee{}{X_n(t \wedge \tau^r) \Rightarrow X(t \wedge \tau^r)}
\n for all $r$. Also $r \rightarrow \infty$ implies $\tau^r \rightarrow \infty$. Therefore $X_n \rightarrow X$.
}

%*************************************************
%*************************************************
\section{Conclusion}
%*************************************************
%*************************************************
We have derived the properties for the operator representation of the multivariate dependence theorem attributed to Sklar which describes the unique representation of general dependence structures linking marginal distributions. This was developed in the context of multidimensional correlated Markov diffusion processes, in the process we defined the copula infinitesimal generator which we interpret as a functional copula specification of Sklars theorem extended to representations involving multivariate generalized diffusion processes. This allowed us to develop a copula function mapping framework which we demonstrated can be accurately and efficiently obtained under a discretization scheme proposed. 

We als demonstrated that it is possible to represent, in both continuous and discrete space, that a multidimensional correlated generalized diffusion is a linear combination of processes that originate from the decomposition of the starting multidimensional semimartingale. This result not only reconciles with the existing theory of diffusion approximations and decompositions, but defines the general representation of infinitesimal generators for both multidimensional generalized diffusions and as was demonstrated allows for a new functional copula characterization of copula density dependence structures. This new result provides immediate representation of the approximate solution for correlated stochastic differential equations. We demonstrate desirable convergence results for the proposed multidimensional semimartingales decomposition approximations in both a strong sense and weak sense, including explicit expressions for the rate of convergence of the discretized approximations.\\

\end{document}